\documentclass[opre,final,nonblindrev]{informs3_for_arxiv}

\pdfoutput=1

\OneAndAHalfSpacedXI 
\usepackage{hyperref}

\newcommand{\III}{\mathcal{I}}
\newcommand{\JJJ}{\mathcal{J}}
\newcommand{\KKK}{\mathcal{K}}

\newcommand{\OOO}{\mathcal{O}}

\usepackage{endnotes}
\let\footnote=\endnote
\let\enotesize=\normalsize
\def\notesname{Endnotes}%
\def\makeenmark{\hbox to1.275em{\theenmark.\enskip\hss}}
\def\enoteformat{\rightskip0pt\leftskip0pt\parindent=1.275em
  \leavevmode\llap{\makeenmark}}
\usepackage{algorithm}
\usepackage{algpseudocode}
\usepackage{multirow}
\usepackage{subfig}
\usepackage{graphicx}
\usepackage{bm}

\usepackage{natbib}
 \bibpunct[, ]{(}{)}{,}{a}{}{,}%
 \def\bibfont{\small}%
\TheoremsNumberedThrough     
\ECRepeatTheorems

\EquationsNumberedThrough    

\begin{document}




\TITLE{Online Resource Allocation for Reusable Resources}

\ARTICLEAUTHORS{%
\AUTHOR{Xilin Zhang}
\AFF{National University of Singapore, \EMAIL{zhangxilin@u.nus.edu}}
\AUTHOR{Cheung Wang Chi}
\AFF{National University of Singapore, \EMAIL{isecwc@nus.edu.sg}} 
} 

\ABSTRACT{We study a general model on reusable resource allocation under model uncertainty. A heterogeneous population of customers arrive at the decision maker's (DM's) platform sequentially. Upon observing a customer's type, the DM selects an allocation decision, which leads to rewards earned and resources occupied. Each resource unit is occupied for a random duration, and the unit is available for another allocation after the usage duration. Our model captures numerous applications involving admission control and assortment planning. The DM aims to simultaneously maximize multiple types of rewards, while satisfying the resource constraints and being uncertain about the customers' arrival process. We develop a near-optimal algorithm that achieves $(1-\epsilon)$ fraction of the optimal expected rewards, where the error parameter $\epsilon$ decays to zero as the resource capacity units and the length of the horizon grow. The algorithm iteratively applies the Multiplicative Weight Update algorithm in a novel manner, which balances the trade-off among the amounts of rewards earned, resources occupied and usage durations.  

}



\maketitle

%


\section{Introduction}

Resource allocation problems constitute a central class of problems in Operations Research, with diverse applications in supply chain management and revenue management. In a resource allocation problem instance, the decision maker (DM) assigns limited inventories of resources to a population of customers, who arrive sequentially during the planning horizon. When a customer arrives, the DM makes an assignment decision. The assignment decision is a generic term, which has different forms in different applications. Examples include admission decisions, pricing decisions, matching decisions, and assortment planning decisions. The DM's decision together with the customer's reaction lead to a set of allocation outcomes, consisting of the amounts of rewards earned and the amounts of resources consumed. The DM aims to maximize the total rewards, while satisfying the capacity constraints.

In the rise of the sharing economy, the allocation of reusable resources has emerged to be an important research topic for contemporary applications. A resource is said to be reusable, if a unit of resource, once assigned, becomes unavailable for a period of usage duration. Upon completion of usage, that unit can be assigned to another customer. Unlike the traditional setting, the set of allocation outcomes consists of the amount of reward earned, the amounts of resources consumed, as well as the usage durations of the consumed resource units. Examples of real-life reusable resource allocation problems include hotel booking, car renting, cloud computing services and emergency department management. Therefore, it is worthwhile to study reusable resource allocation problems considering their wide range of real-life applications.

We incorporate the following four features in our reusable resource allocation model. 
\begin{enumerate}
\item[(F1)] \textit{Multiple objectives.} The DM aims to simultaneously maximize multiple kinds of rewards such as the revenue, the market share and the service level. 
\item[(F2)] \textit{Customer heterogeneity.} Given the same allocation decision, different types of customers exhibit different allocation outcomes. Each allocation outcome is characterized by a stochastic vector consisting of the resources consumed, the rewards earned, and the usage duration.
\item[(F3)] \textit{Online setting.} \textcolor{black}{During the planning horizon, the DM makes allocation decision contingent upon each arriving customer's type, as well as the mean allocation outcome of that customer type. }
In contrast, the DM does \textit{not} know \textcolor{black}{the customer arrival process, which is characterized by the probability distribution of each arrival's type.}
\item[(F4)] \textit{Finite usage duration.} An allocation decision leads to reduction of available resource units over a stochastic usage duration,
instead of permanently reducing the available resources in the remaining planning horizon.
\end{enumerate}
The features (F1-F3) are shared by both non-reusable and reusable resource allocation problems. In contrast, the last feature (F4) of finite usage duration is a distinct feature of the reusable resource setting. As detailed in our forthcoming literature review in Section \ref{sec:lit}, a stream of research works on non-reusable resource allocation has been conducted with an increasing level of generality, culminating to the works of \citep{devanur2019near,balseiro2022best} that achieve near-optimality under general models that incorporate features (F1-F3) in non-reusable resource allocation settings. In this paper, we answer the following research question affirmatively in the related context of reusable setting:
\vspace{0.25cm}
\begin{center}
    \textbf{Question}: Can we achieve near-optimality in a general model that incorporates \\
     features (F1 - F3), in the context of \emph{reusable} resource allocation (F4)?
\end{center}
\vspace{0.25cm}
While the reusable setting is a generalization of the non-resuable setting, the algorithm design and analysis for the latter does not carry over straightforwardly to the former. Indeed, in the reusable setting, the amounts of in-use resource units fluctuate over the entire planning horizon, while the algorithm design and analysis in the non-reusable setting rely heavily on on the monotonically decreasing inventories. In addition, as surveyed in our subsequent literature review, the development on online reusable resource allocation is still in its infancy as compared to the case of online non-reusable resource allocation \citep{devanur2019near,balseiro2022best} . In particular, this paper is the first to study the general settings \citep{devanur2019near,balseiro2022best}  in the context of reusable resources, hence bridging a major gap between the research works on non-reusable and reusable resources.

To bridge the gap, this paper proposes to study an online reusable resource allocation model that incorporates all the four features (F1-F4). Different types of customers arrive sequentially in discrete time steps, and the sequence of customer types is governed by a time homogeneous multinomial process. The length of the planning horizon is not known. The type of customer arrival in each time step follows a common, but latent, probability distribution. Our goal is to find a policy for deciding on the action to take in each 
time step, so that the total rewards earned over the planning horizon are maximized. The main novel contributions of our work are elaborated as follows.
\begin{itemize}
\item \textit{Model generality.} We proposes a general reusable resource allocation model, where the rewards, the customer types, and the DM's allocation actions can be customized to model a variety of decisions in various applications, including admission control, matching, pricing, and assortment planning. Our formulation allows multi-objective optimizations, where the agent aims to minimize the gap between his vector of objective values and a pre-specified vectorial target. In addition, we allow at any time step, the random variables of the amounts of rewards generated, the amounts of resources consumed, and the usage duration to be arbitrarily correlated. This generalize the resource consumption model in related research works on assortment planning, which either one unit of a resource or no resource is consumed in each time step. We highlight that the correlation between resource consumption amounts and usage durations is essential in modeling workforce management applications, where assigning more working power usually leads to a shorter completion time.
\item \textit{Algorithm design.} We develop an iterative multiplicative weight update (iMWU) algorithm. In each phase, given the historical observations from past phases, we run a MWU algorithm which produces weight vectors that trade-off among the amounts of rewards earned, the amounts of resources consumed as well as the usage durations. From these weight vectors, we sample one vector for each customer arriving in the current phase that determines the allocation action. We quantify the performance of our policy by comparing our rewards with a suitably chosen offline linear program benchmark that upper bounds the optimal dynamic policy.
\item \textit{Near-optimal algorithm performance.} For any given confidence parameter $\delta \in (0, 1)$, we show that with probability at least $1-\delta$, our online algorithm achieves a time-averaged reward at least $1- (3\sqrt{\xi\log(|\III_c|/\xi)} + O(\frac{\log(T)}{\sqrt{T}}))$ fraction of an offline benchmark that upper bounds the expected optimum. In the approximation ratio, $\xi \in (0,1)$ is the largest ratio of the amount of resource consumed in a time step to the total amount of the resource, $|\III_C|$ is the number of resources, and $T$ (unknown to the DM) is the length of the planning horizon. The $O(\cdot)$ notation hides a multiplicative factor that does not depend on $T$, but only depends on other model parameters and $\delta$. We highlight that a loss factor of $O(\sqrt{\xi})$ loss factor is inevitable, since \cite{agrawal2014fast} show that it is unavoidable in the non-reusable resource setting, a special case of ours when the usage durations for all customers are $T$ with certainty. Remarkably, the approximation factor does not deteriorate as the number of customer types or the number of actions increases. The former property makes our algorithmic framework attractive in applications involving many customer types and complex actions. 

One of the major challenge in our algorithm design and analysis is the \emph{structural non-stationarity} on the resource consumption model, where the amounts of resource occupied in different time steps are differently distributed, due to the heterogeneity induced by random usage durations. The non-stationarity makes existing algorithms for non-reusable resouce allocation inapplicable, and motivates our alternative algorithmic approach.
\item \textit{Numerical validation.} On top of the theoretical analysis, we conduct numerical experiments on an assortment planning problem where both the number of customer types and the number of DM’s actions are large. A column-generation-based decomposition approach is used to facilitate computation efficiency of our policy. The results exhibit the strong numerical performance of our algorithm. 
\end{itemize}
To this end, we highlight that our work is the first to develop near-optimal algorithms for online stochastic reusable resource allocation problems, where the usage durations are governed by arbitrary probability distributions.

\section{Related literature}\label{sec:lit}
We provide an overview of existing research works on resource allocation problems with non-reusable and reusable resources, in both the offline and the online settings. The problem contexts include admission control, pricing, assortment planning and queueing systems. For resource allocation problems, the offline setting and the online setting differs in whether there exists model uncertainty. In an offline (or Bayesian) setting, there is no model uncertainty. Although the DM does not know the exact realization of the customer arrival process, he does know all the underlying distributions a priori. On the other hand, in an online setting, the DM does not know the full model on the customer arrival process. He has no knowledge of the distributions of certain parameters beforehand. These parameters can be either stochastic or adversarial depending on the nature of the environment. In an online stochastic scenario, the parameters are (often i.i.d.) drawn from some unknown distribution. In an online adversarial scenario, the parameters can be seen as chosen by an adversary. 

\subsection{Non-reusable resource allocation problems}
Traditional non-reusable resource allocation problems have been extensively studied in both the offline setting \citep{adelman2007dynamic, alaei2012online,Jasin15,bumpensanti2020re, vera2021bayesian} and the online setting \citep{FeldmanMMM09, FeldmanHKMS10, agrawal2014fast,Jasin14, devanur2019near, LiY19,balseiro2020dual}. In offline settings, recent progress on dynamically adjusting the decisions based on the inventory level leads to $O(1)$ regret bounds in \citep{Jasin15,bumpensanti2020re, vera2021bayesian}. In online stochastic settings, \citet{AgrawalWY14} show that no algorithm can achieve a better performance than a $1 - O(\sqrt{\xi})$ fraction of the optimum. A cornucopia of online stochastic algorithms are proposed, and these algorithms involve adaptive weighing processes that trade-off between the rewards earned and the resources consumed. For example, \citet{agrawal2014fast} propose a dual-based algorithm for online linear optimization, and demonstrate the time average approachability of constraints. \citet{balseiro2020dual} develop a primal-dual learning algorithm based on analyzing the stopping time when the inventories are depleted. \citet{devanur2019near} construct an online policy for a general multi-objective resource allocation model. 
All these three papers achieve a performance guarantee of at least $1 - \tilde{O}(\sqrt{\xi})$ times their respective optima, where $\tilde{O}(\cdot)$ hides a $\log$ multiplicative factor in terms of $\xi$ or the number of resources. In our forthcoming main result statements, we demonstrate that our proposed algorithm also achieves a factor of $1 - \tilde{O}(\sqrt{\xi})$ times the optimum, when $T$ is large. 

Despite the abundance of researches on the non-reusable resource allocation problems, these works depend heavily on the property that the amounts of available resources are non-increasing across time. Therefore, they cannot be directly generalized to reusable resource allocation models, where the amounts of available resources can fluctuate across time. On one hand, the amounts of available resources decrease when they are assigned to a customer. On the other hand, the amounts of available resources increase when the in-use units are returned. \cite{devanur2019near} provide an algorithm on resource allocation with model non-stationarity, which could potentially handle reusable resource settings Nevertheless, they require the crucial prior knowledge on the expected number of in-use units under the optimal dynamic policy, which could be often unavailable in real-life applications. The reason why they need extra information can be attributed to the fact that with fluctuating inventories, the DM has to keep track of when each resource unit is allocated and returned. A crucial contribution of our work is on showing that the structural property entailed by the reusable resource allocation model can be harnessed to bypass the need of the previous mentioned prior information, albeit with a significantly different algorithm framework. 
 
Lastly, it is also worth mentioning that for non-reusable resource allocation in the online setting, the DM must know the horizon $T$ a priori to obtain near-optimal performance. In contrast, we show that in our online problem concerning reusable resources, $T$ needs not be known to the decision maker.

\subsection{Reusable resource allocation problems in the offline setting}
Reusable resource allocation problems are less well understood compared to traditional (non-reusable) resource allocation problems. Most works concerning reusable resources focus on the offline setting. \citet{levi2010provably} consider an admission control problem and a pricing problem with a single reusable resource. 
They propose an LP guided policy, and establish an  approximation ratio $1 - O(\sqrt{\xi})$ for their policy in an average reward and infinite horizon setting. Specifically, their approximation ratio approaches $1 - \sqrt{2\xi/\pi}$ when $\xi$ tends to 0. \citet{chen2017revenue} study a similar setting with advance reservation. They consider a scaling regime, which requires synchronously scaling the system capacity and the customer arrival rate at the linear rate of $n$. They show that the optimality gap shrinks when $n$ grows to infinity. \citet{lei2020real} study a general network pricing model with non-stationary customer arrivals and advance reservations in a similar scaling regime, where the usage durations of resources are deterministic. 
Compared with these works, our work achieves an approximation ratio of $1- \tilde{O}(\sqrt{\xi})$ in general models with large capacities, and our result does not require the scaling regime assumption as in \cite{chen2017revenue,lei2020real}.

While some offline papers regarding reusable resources focus on achieving near-optimal performance in large-scale systems, some other works derive universal results independent of the system scale. \citet{rusmevichientong2020dynamic} study an assortment planning problem based on a dynamic programming formulation. 
 They design a policy based on affine approximation, and the policy accrues at least $\frac{1}{2}$ of the expected revenue compared with the optimum of the underlying dynamic programming problem. \cite{baek2022bifurcating} study a network revenue management problem with reusable resources and derive a policy with approximation ratio that depends on the sizes of requests.  
 \citet{feng2020near} design a class of discarding policies on an assortment planning problem with reusable resources. Their algorithm achieve an approximation ratio of $1 - \min\left\{\frac{1}{2}, \sqrt{\xi \log \frac{1}{\xi}}\right\}$. \citet{besbes2021static} study a multi-objective pricing problem and prove that a well-chosen static pricing policy guarantees $78.9\%$ of the optimum. 

Another line of offline works study admission control or pricing problems in a queueing system where the reusable resources are regarded as servers, and the customers are regarded as jobs. The jobs arrive sequentially according to a stationary poisson process. They usually consider the a similar scaling regime (termed heavy-traffic regime in queueing literature) as in \cite{chen2017revenue, lei2020real} where the server capacity and the job arrival rate scale at rate $n$. Correspondingly, they provide performance guarantees when the scaling parameter $n$ grows. Most queueing works adopt a static policy guided by the fluid approximation of the system. \citet{kim2018value} uses a two-price policy based on different system states and achieves a better performance guarantee (faster converge rate to optimality regarding $n$). It is worth mentioning that in queueing systems, jobs can often wait before being served. Our problem focus on the loss system, which is a special case in queueing systems where the customers are lost immediately if there is no idle server. 

\subsection{Reusable resource allocation problems in the online setting}
For the online setting concerning reusable resource allocation problems, most works focus on the adversarial scenario. The performance of an algorithm is usually quantified by its competitive ratio. An online algorihtm is said to be $\alpha$-competitve, if it achieves an expected reward of at least $\alpha$ times the optimum on any instance. In this line of works, the customer arrival process is allowed to be arbitrary, but the usage durations are assumed to be i.i.d. and independent of the customer types, in order for achieving non-trivial competitive ratios. \citet{gong2019online} study an assortment planning problem, and show that the myopic policy is $(1/2)$-competitive. For each customer, they offer the assortment with the largest expected revenue among all assortments consisting of available resources, and their algorithm is oblivious to the probability distributions of the usage durations. 
The research work  \citep{feng2021online} considers a assortment planning model on reusable resources, by incorporating the novel feature of exogenous inventory replenishment, and they achieve a constant competitve ratio by designing online algorithm based on inventory balancing. \citet{goyal2020asymptotically} and \citet{goyal2020online} develop fluid approximation guided algorithms, which achieve a $1-\frac{1}{e}$ competitive ratio for the case random usage duration. Usually, the online adversarial works do not require knowing the length of the planning horizon $T$ a priori. Compared with these works, the rewards of our algorithm converge to the optimum when $\xi$ decreases under the assumption of stochastic online arrivals, in contrast to the adversarial settings where the competitive ratio remains constant. 

In the stochastic scenario, \cite{kanoria2019blind} study online resource allocation problems in closed networks, where the resource units circulate among different nodes in a networked system. They propose a novel variant of the backpressure policy, where there algorithm's performance converges to the optimum when the number of resource units and the time horizon grow. Crucially, they assume that one resource unit is relocated at a time, thus in their model the usage duration is deterministic with value 1. While our model does not consider a general networked system as in \citep{kanoria2019blind}, our model allows allocating multiple resources to each customer and general usage durations.

\textbf{Notation. }We denote the set of non-negative real numbers as $\mathbb{R}_{\geq 0}$, and the set of strictly positive real numbers as $\mathbb{R}_{>0}$. Further, we denote $\mathbb{R}^n_{\geq 0} = \left\{x \in \mathbb{R}^n : x \geq 0 \right\}$. Occasionally, we use $[T]$ as a shorthand of $\{1, \ldots, T\}$. For an optimization problem (Z), we denote opt(Z) as its optimum. A table of important notation for our model and our algorithms is provided in Appendix section \ref{tab:notation}.

\section{Problem formulation}
In this section, we introduce a discrete-time online stochastic model for allocating reusable resources. Our model incorporates multi-objective optimization and customer heterogeneity, and accommodates a variety of allocation decisions such as admission control and assortment planning decisions.  

\subsection{Model}
\textbf{Rewards and resources. }The decision maker (DM) aims to simultaneously maximize multiple types of rewards, by allocating multiple types of reusable resources to a population of customers. The reward types and the resource types are respectively indexed by two disjoint finite sets $\III_r$ and $\III_c $. A generic reward type or resource type is denoted as $i$. For each $i\in \III_c$, the DM has $c_i\in \mathbb{R}_{>0}$ units of resource $i$ for allocation. Each customer is associated with a customer type $j\in \JJJ$, which reflects the customer's characteristics, such as his/her gender, age, interest, etc.  
We denote the set of all possible allocation decisions as $\KKK$. For brevity, we call $\KKK$ the action set, and each element $k\in \KKK$ as an action. The action set can be used to model a broad range of decisions, such as the product assortment to offer. We elaborate on these decisions in Section \ref{sec:app}. 


The DM allocates the resources in $T$ discrete time steps. In time step $t \in \{1,\ldots,T\}$, at most one customer arrives. We denote the customer type of the arrival at time $t$ as $j(t)$. In particular, we designate the type $j_{\textsf{null}}$ (which is assumed to be in $\JJJ$) to represent the case of no arrival. 
We assume that $j(1), \ldots, j(T)$ are independently and identically distributed (iid) random variables over $\JJJ$. We denote $p_j = \Pr(j(1) = j)$, and $\mathbf{p} = \{p_j\}_{j\in \JJJ}$. 

When a customer (denote his type as $j$) arrives, the DM chooses an action $k\in \KKK$. The choice leads to an array of stochastic outcomes $(W_{jk}, A_{jk}, D_{jk})\in \mathbb{R}^{\III_r}_{\geq 0} \times \mathbb{R}^{\III_c}_{\geq 0} \times \mathbb{Z}^{\III_c}_{\geq 0} $. In $W_{jk} = (W_{i jk})_{i \in \III_r}$, the random variable $W_{i j k}$ is the amount of type-$i$ reward earned, for each $i\in \III_r$. In $A_{jk} = (A_{ijk})_{i\in \III_c}$, the random variable $A_{ijk}$ is the amount of type-$i$ resources occupied, for each $i\in \III_c$. In $D_{jk} = (D_{ijk})_{i \in \III_c}$, the random variable $D_{ijk}$ is the usage duration of the $A_{ijk}$ type-$i$ resource units, for each $i\in \III_c$. 



For the no arrival customer type $j_{\textsf{null}}$, we stipulate that $\Pr(W_{i', j_{\textsf{null}}, k} = A_{i, j_{\textsf{null}}, k}= D_{i, j_{\textsf{null}}, k} = 0) = 1$ for all $i' \in \III_r, i\in \III_c, k\in \KKK$, since there should be no reward earned and no resource occupied in the case of no arrival. To ensure feasibility in our resource constrained model, we assume that there exists a null action $k_{\textsf{null}}\in \KKK$ that satisfies $\Pr(W_{i', j, k_{\textsf{null}}} = A_{i, j, k_{\textsf{null}}}= D_{i, j, k_{\textsf{null}}} = 0) = 1$ for all $i' \in \III_r, i\in \III_c, j\in \JJJ$. Selecting the null action is equivalent to rejecting a customer, whereby no reward is earned and no resource unit is occupied. 

We denote $\OOO_{jk}$ as the joint probability distribution of $(W_{jk}, A_{jk}, D_{jk})$, with shorthand $(W_{jk}, A_{jk}, D_{jk}) \sim \OOO_{jk}$. We allow the entries $W_{1jk}, \ldots, W_{|\III_r|,jk}$ in the vectorial outcome $W_{jk}$ to be arbitrarily correlated, and similarly for $A_{jk}, D_{jk}$. In addition, we allow the random vectors $W_{jk}, A_{jk}, D_{jk}$ to be arbitrarily correlated. 
We assume that $W_{ijk} \in [0,w_{\max}]$ almost surely for each $i \in \III_r, j\in \JJJ, k\in \KKK$, and we assume $A_{ijk} \in [0, a_\text{min}]$, and $D_{ijk} \in \{0,1, \ldots, d_{\max}\}$ almost surely for each $i \in \III_c, j\in \JJJ, k\in \KKK$. 

Additionally, we denote $w_{i jk} = \mathbb{E}[W_{i jk}]$, $a_{ijk} = \mathbb{E}[A_{ijk}]$, and $d_{ijk} = \mathbb{E}[D_{ijk}]$, and denote $w_{jk} = (w_{i jk})_{i\in \III_r}, a_{jk} = (a_{ijk})_{i\in \III_c}, d_{jk} = (d_{ijk})_{i\in \III_c}$. Crucially, we define
$$
v_{ijk} = \mathbb{E}[A_{ijk}D_{ijk}],
$$
which can be interpreted as the expected volume of resource $i$ consumed under customer type $j$ and action $k$. For example, in the case of cloud computing, the quantity $v_{ijk}$ could represent the expected type $i$ resource-hour (for example, CPU-hour when $i$ represents CPUs) needed for a type-$j$ computational task under allocation decision $k$. In the case of human resource management in healthcare settings, the quantity $v_{ijk}$ could represent the expected type $i$ man-hour (for example, nurse-hour when $i$ represents nurses) for a type $j$ patient under allocation decision $k$. Since the random variables $A_{ijk}, D_{ijk}$ can be correlated, in general $v_{ijk}$ needs not be equal to $a_{ijk}d_{ijk}$. We denote $v_\text{max} = \max_{i,j,k}\{v_{ijk}\}$.


\textbf{Dynamics and non-anticipatory policies. }At each time step $t\in \{1, \ldots, T\}$, three events happen. Firstly, the DM observes the type $j(t)\sim \textbf{p}$ of the time $t$ customer, and the mean outcomes $\{(w_{j(t), k}, v_{j(t), k})\}_{k\in \KKK}$ specific to the type $j(t)$. Secondly, the DM chooses an action $k(t)\in \KKK$.  Thirdly, the DM observes the array of stochastic outcome $(W(t), A(t), D(t))\sim \mathcal{O}_{j(t), k(t)}$. Altogether, the observation at time step $t$ consists of $j(t), \{(w_{j(t), k}, v_{j(t), k})\}_{k\in \KKK} , W(t), A(t), D(t)$. We denote the entries in the stochastic vectorial outcomes as $W(t) = (W_i(t))_{i\in \III_r}, A(t) = (A_i(t))_{i\in \III_c}, D(t) = (D_i(t))_{i\in \III_c}$.  

For each $i\in \III_r$, the DM earns $W_i(t)$ units of type $i$ reward. For each $i\in \III_c$, the time $t$ customer occupies $A_{i}(t)$ units of resource $i$ during time steps $t, t+1, \ldots, t + D_i(t) -1$.  These $A_{i}(t)$ units of resource $i$ would be available for another allocation again from time step $t+ D_i(t)$ onwards. 

The DM chooses the actions $k(1), \ldots, k(T)$ by a non-anticipatory policy. The choice of $k(t)$ is only based on (a) the observed type $j(t)$, (b) the historical observations $\mathcal{H}(t-1)$ from time 1 to time $t-1$, (c) the internal randomness of the DM. By contrast, the action $k(t)$ is chosen without knowing the future customer arrivals $j(t+1), \ldots, j(T)$.

\textbf{Objective. }The DM aims to
maximize $\mathbb{E}\left[\min_{i \in \III_r}\{ \sum^T_{t=1} W_i(t)  \}/T\right]$, subject to resource constraints and model uncertainty. The quantity $\sum^T_{t=1} W_i(t) /T$ is the average type $i$ reward. The maximization objective achieves the simultaneous maximization 
of all the reward types by ensuring  max-min fairness. For each resource $i\in \III_c$ and each time step $t\in \{1, \ldots, T\}$, we require that the resource constraint
\begin{equation}\label{eq:constraint}
\sum^t_{\tau = 1}  \mathbf{1}( D_{i}(\tau) \geq t - \tau + 1) A_{i}(\tau) \leq c_i
\end{equation}
 holds with certainty. The left hand side in (\ref{eq:constraint}) represents the amount of occupied type $i$ resources at time step $t$. In particular, the time $\tau$ customer occupies $\mathbf{1}( D_{i}(\tau) \geq t - \tau + 1) A_{i}(\tau)$ units of type $i$ resource at time step $t$.

Finally, by specializing $\Pr( D_{ijk}=T)=1$ for all $i\in \III_c, j\in \JJJ\setminus \{j_{\textsf{null}}\}, k\in \KKK\setminus \{k_{\textsf{null}}\}$, our model specializes to the non-reusable resource allocation model in online stochastic settings studied in \citep{devanur2019near}, where an allocated resource unit cannot be re-allocated again in the planning horizon. Such settings are studied in a cornucopia of research works \citep{GoelM08,AgrawalWY14,agrawal2014fast,devanur2019near,li2021online, balseiro2022best}. 

\textbf{Model uncertainty.} The DM is uncertain about the model, in the sense that he does not know the probability distribution $\mathbf{p}$ over the customer types, and he does not know the horizon length $T$. While the collection $\{(w_{j(t), k}, v_{j(t), k})\}_{k\in \KKK}$ of the mean outcomes specified to the time $t$ customer is only revealed sequentially, the bounds $a_\text{max}, d_\text{max}, w_\text{max}, v_\text{max}$ are known to the DM before the online process begins. 

Our model uncertainty scenario includes the case when the DM knows the mean outcomes $\{(a_{jk}, d_{jk}, w_{jk}, v_{jk})\}_{j\in \JJJ,k\in  \KKK}$ before the online dynamics begin, but does not know $\textbf{p}, T$. Such an uncertainty scenario represents the case when the DM knows the relationship between a customer's type and his/her preferences, but the DM does not know the number of customers belonging to each customer type. For example, the DM could have acquired the knowledge on the type-preference relationship through his previous interaction with another customer population. While the DM could carry over his knowledge on $w_{i jk}, a_{ijk}, d_{ijk}, v_{ijk}$ when he faces a new customer population, the DM is still required to overcome the model uncertainty on $\textbf{p}$, since the new population's composition could differ from the previous population. We elaborate more on this in Section \ref{sec:app}.
 
\textbf{An Online Integer Program Interpretation.} To facilitate our technical discussions, we rephrase the online resource allocation problem as an online integer program. For a non-anticipatory policy $\pi$ that is used to select the actions $k(1), \ldots, k(T)$, we let the binary decision variable $X^\pi_k(t)\in \{0, 1\}$ be the indicator variable on the event that action $k$ is chosen at time $t$. That is, $X^\pi_k(t) = 1$ if and only if $k = k(t)$. In addition, for each $j\in \JJJ, k\in \KKK$, let $(W_{j k}(t), A_{j k}(t), D_{jk}(t))^T_{t=1}$ be $T$ iid samples drawn from $\mathcal{O}_{jk}$, and recall our definition that $j(1), \ldots, j(T)$ are iid with the common probability distribution $\textbf{p}$. The DM's objective can be equivalently phrased as the following online stochastic integer program, where all constraints are to be satisfied with certainty:
\begin{subequations}
\begin{alignat}{2}
\text{(IP-C) }\max\limits_{\text{non-anticipatory }\pi}  &~\mathbb{E}[\lambda^C] & \nonumber\\
\text{s.t.}  & ~\frac{1}{T}\sum^T_{t=1} \sum_{k\in \KKK} W_{i, j(t), k}(t) X^\pi_{k}(t) \geq  \lambda^C     &\quad &\forall i\in \III_r   \nonumber\\
&\sum^t_{\tau = 1} \sum_{k\in \KKK} \mathbf{1}( D_{i,j(\tau),k}(\tau) \geq t - \tau + 1) A_{i,j(\tau),k}(\tau) X^\pi_{k}(\tau) \leq c_i       &\quad & \forall i\in \III_c,~ t\in [T] \nonumber\\
& \sum_{k\in \KKK} X^\pi_k(t) = 1 &\quad &\forall ~ t\in [T] \nonumber\\
&X^\pi_k(t) \in \{0,1\}      &\quad &\forall k\in \KKK,~ t\in [T] \nonumber.
\end{alignat}
\end{subequations}
By requiring $\pi$ to be non-anticipatory, we are requiring that $X^\pi_{k}(t)$ is $\sigma(\{j(t)\} \cup {\cal H}(t-1)\cup \{U(t)\})$-measurable, where the random variable $U(t)$ represents the internal randomness of the DM used at time $t$. It is worth mentioning that by setting $D_{ijk} = T$ with certainty for all $i,j,k$, we recover the non-reusable resource allocation model in \citet{devanur2019near}, which is a multi-objective generalization of \citet{AgrawalWY14} and \citet{li2021online}. Although our model only change one random variable compared with the non-reusable setting, the structure of the model changes entirely. We have $T|\III_c|$ resource constraints to keep track of the available resources in each time step $t \in [T]$. Moreover, the resource constraints involve non-stationariety on resource consumption, even though we assume a time homogeneous customer arrival process.

We further discuss on the non-stationarity on resource consumption. 
Suppose a type-1 customer would occupy one unit of allocated resource for a random duration of $D$ time units, where $\Pr(D = 5)= 0.3, \Pr(D =10) = 0.7$, upon allocation decision $k$ (We omit the subscripts to ease the notation in the exposition). Consider the case when one type-1 customer arrives during time steps 1 and 6, and the decision maker selects decision $k$ in both time steps 1 and 6. Let's consider the amount of resources units in use by these two customers at time step 6. On one hand, the stochastic outcomes to these two customers are identically distributed, since the customer type and the allocation decision are the same in the two time steps. On the other hand, the amount of resource units in use by the time-1 customer is $\sim \text{Bern}(0.7)$, while the amount of resource in use by the time-6 customer is equal to 1 with certainty. This intrinsic structural non-stationarity, and the fluctuating amounts of in-use resource units, 
prevent direct adaptation of the existing algorithms on non-reusable algorithms. The above-mentioned technical challenge motivates our alternate algorithm design and analysis in our subsequent discussions. 

\textbf{Our Goal. }The optimal value  of (IP-C) is the optimal expected reward obtained by any non-anticipatory algorithm. However, the problem (IP-C) is intractable in general, due to the curse of dimensionality. Thus, our aim is to achieve near-optimality. More precisely, for any given confidence parameter $\delta\in (0, 1)$, we seek to construct a policy $\bar{\pi}$ that achieves
$$
\frac{1}{T}\sum^T_{t=1} \sum_{k\in \KKK} W_{i, j(t), k}(t) X^{\bar{\pi}}_{k}(t)\geq  \alpha \cdot \text{opt(IP-C)} - O\left(\frac{1}{T^\beta}\right) \qquad \text{for each $i\in \III_r$},
$$ 
with probability at least $1-\delta$, for some parameters $\alpha, \beta \in (0, 1]$. The parameter $\alpha$ is an approximation factor. In the case of single objective (i.e. $|\III_r| = 1$) and full model certainty (meaning that $\textbf{p}, T$ are known), existing works have established approximation ratios of $1 - \min\left\{\frac{1}{2}, \sqrt{\xi \log \frac{1}{\xi}}\right\}$ \citep{feng2020near}, where $\xi = a_\text{max} / c_\text{min}$. More precisely, \cite{feng2020near} consider the case of maximizing the expected reward instead of the high probability objective, and they set $a_\text{max} = 1$. On the specialization of our model to non-reusable resources ($D_{ijk} = T$ with certainty for all $i, j, k$), \citep{devanur2019near,agrawal2014fast} also achieve an approximation ratio of $1 - O(\sqrt{\xi})$. Next, the error term $O\left(\frac{1}{T^\beta}\right)$ is due to the randomness of the stochastic outcome and the uncertainty over $\textbf{p}$ and $T$.


While our optimization problem (IP-C) promotes max-min fairness the accrued rewards among different types, (IP-C) can also be adjusted to model the following \emph{KPI optimization problem}. Specifically, the DM is endowed with an array of KPI parameters $\boldsymbol{\sigma} = \{\sigma_i\}_{i \in \III_r}$, and he aims to achieve 
$$
\frac{1}{T} \sum^T_{t=1} \sum_{k\in \KKK} W'_{i, j(t), k}(t) X^{\bar{\pi}}_{k}(t) \geq  \alpha \cdot \sigma_i \cdot \text{opt(IP-C)} - O\left(\frac{1}{T^\beta}\right) \qquad \text{for each $i\in \III_r$}
$$
with probability at least $1 - \delta$, for some random reward variables $W'_{ijk}$. In this case, by letting $W_{i, j(t), k}(t) = W'_{i, j(t), k}(t) / \sigma_i$, we incorporate the KPI into our model, and recover the exact same formulation and goal as (IP-C).
 


\textbf{A Tractable Upper Bound of (IP-C).} 
While it is desirable to construct a non-anticipatory algorithm that achieves (or nearly achieves) the optimal value of (IP-C), the optimal value of (IP-C) is analytically intractable due to the curse of dimensionality. The intractability motivates us to consider an alternative LP, dubbed (LP-E). We provide (LP-E) in the following, where the realization of the customer arrivals, their usage duration and outcomes exactly follow the expectation:
\begin{subequations}
\begin{alignat}{2}
\text{(LP-E) }\max  &~\lambda^E & \nonumber\\
\text{s.t.}  &\sum^T_{t=1} \sum_{j \in \JJJ} \sum_{k\in \KKK} p_j w_{ijk} x_{jk}(t)\geq T \lambda^E     &\quad &\forall i\in \III_r   \nonumber\\
&\sum^t_{\tau = 1}  \sum_{j\in \JJJ}\sum_{k\in \KKK} p_j \mathbb{E}[A_{ijk} \mathbf{1}(D_{ijk} \geq t - \tau + 1)] x_{jk}(\tau) \leq c_i       &\quad & \forall i\in \III_c,~ t\in [T] \nonumber\\
&\sum_{k\in \KKK} x_{jk}(t)\leq 1      &\quad &\forall j\in \JJJ,~ t\in [T] \nonumber\\
&x_{jk}(t)\geq 0      &\quad &\forall j\in \JJJ, ~k\in \KKK,~ t\in [T] \nonumber.
\end{alignat}
\end{subequations}
Define the optimal objective value of (LP-E) to be $\lambda^E_*$, and let the optimal objective of (IP-C) be $\lambda^C_*$. We next show the following Lemma.
\begin{lemma}\label{lem:benchmark1}
$\lambda^E_* \geq \lambda^C_*$.
\end{lemma}
Lemma \ref{lem:benchmark1} is proved in Appendix Section \ref{pf:lemma_benchmark1}. 

\subsection{Applications}\label{sec:app}
Our model captures reusable resource allocation in many general settings. We illustrate the relevance of our model to applications through two examples: admission control problems, and assortment planning for reusable resources. 

\textbf{Admission Control. }The DM's action set $\KKK$ consists of two actions, $k_{\textsf{accept}}, k_{\text{null}}$. At the start of each time step, the DM observes a customer's type $j$, associated with expected rewards $\{w_{i,j}\}_{i\in \III_r}$ and expected resource consumption $\{v_{i,j}\}_{i\in \III_c}$. Then, the DM decide to either accept the customer by taking the action $k_{\textsf{accept}}$, or to divert/reject the customer by taking the action $k_\textsf{null}$. In the former case of $k_{\textsf{accept}}$, the DM earns a type $i$ reward $W_{i, j}$ with $\mathbb{E}[W_{i, j}] = w_{i, j}$ (so that $w_{i,j, k_{\textsf{accept}}} = w_{i,j}$). The DM occupies $A_{ij}$ units of type $i$ resource for $D_{ij}$ time steps, with $\mathbb{E}[A_{i, j}D_{i, j}] = v_{i, j}$ (so that $v_{i,j, k_{\textsf{accept}}} = v_{i,j}$). In the latter case of $k_\textsf{null}$, the DM earns no reward for every type, but also does not occupy any resource.

For a more concrete example, consider the case of a DM operating a cloud computing platform during a peak-hour period, which is the planning horizon. The DM could either accept an incoming job, or he could divert the job to a waiting room, where the diverted jobs would only be processed after the peak-hour period. In cloud computing settings, the quantity $v_{i,j,k_{\textsf{accept}}}$ represents the expected number of CPU-hours consumed (let's say $i$ represent the CPU resources) if a type $j$ job is accepted. The same representation generalizes to other resources such as GPUs, RAM, etc. Our general model captures the heterogeneity of jobs through $\JJJ$, where each 
type $j\in \JJJ$ encodes the nature of the job and the computational requirements. In addition, our model captures the need of multi-resource consumption and non-identical resource requirements, through a general model on $\{A_{i,j,k_{\textsf{accept}}}, D_{i,j,k_{\textsf{accept}}}\}_{i\in \III_c}$, where these random outcomes could correlate among each other arbitrarily. Moreover, our multi-objective optimization framework allows the DM to incorporate objectives such as utilization rates of resources and social welfare, in addition to the total revenue earned.

Let us remark on model uncertainty. For computational jobs on training neural networks, an estimate on the amount of computational resources needed could often be obtained based on the job's attributes  \citep{JustusBBM18,ZancatoARBS}. Hence, it is reasonable to assume that $v_{j,\textsf{accept}}, w_{j,\textsf{accept}}$ are known when a type $j$ job arrives. In contrast, given that a type $j$ encodes various attributes of a job, the size of $\JJJ$, which contains all the possible types, could be prohibitively large for the purpose of estimating $p_j$ for each $j\in \JJJ$. Thus, it is desirable to construct an online algorithm that achieves near optimality without knowing $\textbf{p}$.

While the above example is made in the context of cloud computing, our model can also be used in other context such as healthcare resource management and workforce management. In healthcare resource management, resources include hospital beds, doctors and nurses, medical equipment etc. For the emergency department (ED) of major hospitals, incoming patients may be diverted to other hospitals or the outpatient department to mitigate crowding. The quantity $v_{j,k_{\textsf{accept}}}$ represents the expected amounts of healthcare resources occupied over time if a type $j$ patient is admitted into the ED, and the quantity $w_{j,k_{\textsf{accept}}}$ indicates the acuity level of a type $j$ patient ($w_{j,k_{\textsf{accept}}}$ larger for more urgent patients) as well as other objectives of the DM.

\textbf{Assortment Planning }for allocating reusable resources has been studied actively in recent literature \citep{gong2019online,feng2020near,goyal2020online,rusmevichientong2020dynamic}. In assortment planning problems, each resource corresponds to a product for sales. One unit of resource $i$ is associated with a fixed price $r_i$. The DM influences the customers' demands through offering different assortments of resources. Contingent upon the arrival of a customer, say of type $j$, the DM decides the assortment $k \in \KKK$ to display, where $\KKK$ is a collection of subsets of $\III_c$. A popular example for $\KKK$ is the cardinality constrained collection $\KKK = \{k : |k| \leq n , k\subset \III_c\}$, where $n$ represents the maximum number of resources/products that can be displayed by the DM.

Upon being offered an assortment $k$, a customer either chooses a product $i$ in $k$, or chooses no product. We use product 0 to denote the choice of no product. For $i\in k\cup \{0\}$, we let $q_{ijk}$ denote the probability for customer type $j$ to choose product $i$ when offered assortment $k$, and $q_{0,jk}$ denotes the probability of no purchase. Consequently, we have $\sum_{i\in k\cup \{0\}}q_{ijk} = 1$. When a type $j$ customer chooses product $i\in k$, he occupies one unit of resource $i$ for a random duration of $\textsf{Duration}_{ij}$, and none of the not-chosen products is occupied. Thus, we have $A_{ijk} \sim \text{Bern}(q_{ijk})$, and $D_{ijk} = A_{ijk}\cdot \textsf{Duration}_{ij}$. Note that the outcomes $A_{ijk}, D_{ijk}$ are correlated. In the above-mentioned recent works, it is assumed that $A_{ijk},\textsf{Duration}_{ij}$ are independent, so that $d_{ijk} = a_{ijk} \mathbb{E}[\textsf{Duration}_{ij}] = q_{ijk} \mathbb{E}[\textsf{Duration}_{ij}]$.

In the revenue management literature, the choice probability $q_{ijk}$ is often modeled by a random utility choice model. In the following, we provide an illustrative example with the multinomial logit (MNL) choice model, which is one of the most prominent random utility choice models. Each resource $i\in \III_c$ is associated with a feature vector $\boldsymbol{f}_i \in \mathbb{R}^m$, and each customer type $j\in \JJJ$ is associated with a set of feature vectors $\{\boldsymbol{b}_{ij} \}_{i \in \III_c}$, where $\boldsymbol{b}_{ij}  \in \mathbb{R}^m$ for each $i\in \III_c$. The feature vector $\boldsymbol{f}_i$ could involve the fixed price $r_i$, which reflects the intrinsic value of resource $i.$ The random utility of a type $j$ customer on resource $i$ is modeled as $V_{ij} = \boldsymbol{b}_{ij}^\top \boldsymbol{f}_i + \zeta_{ij}$, where $\{\zeta_{ij}\}_{i, j}$ is the collection of iid random variables following the Gumbel distribution. 
 
When an assortment $k$ is offered to a type $j$ customer, the customer either chooses to purchase one unit of product $i_*\in k$, where $i_*\in \text{argmax}_{i\in k} V_{ij}$, when $V_{i^*j} >0$; or, the customer chooses to make no purchase when $V_{i^*j} \leq 0$. A crucial property of the multinomial logit choice model is that
$$
q_{ijk} = \frac{\exp(\boldsymbol{b}_{ij}^\top \boldsymbol{f}_i)}{1+\sum_{\ell \in k} \exp(\boldsymbol{b}_{\ell j}^\top \boldsymbol{f}_\ell)}
$$
if $i\in k$, and $q_{ijk} = 0$ if $i\not \in k$. In complement, the probability of no purchase is  $\frac{1}{1+\sum_{\ell \in k} \exp(\boldsymbol{b}_{\ell j}^\top \boldsymbol{f}_\ell)}$. Consequently, the assortment planning problem with the single objective of maximizing the total revenue \cite{feng2021online,rusmevichientong2020dynamic} can be modeled by setting $\III_r = \{1\}$, with $W_{1, jk} = \sum_{i\in \III_{c}}r_i A_{ijk}$, where $A_{ijk}$ is the Bernoulli random variable with mean $q_{ijk}$. Our model also allows maximizing multiple objective, for example, simultaneously maximizing the revenue of each resource can be modeled by setting $\III_r = \III_c$ and $W_{ijk}=r_i A_{ijk}$ for each $i\in \III_r$. In our forthcoming numerical section, we provide another example of multi-objective optimization where the DM aims to jointly maximize the revenue earned and the sales volumes of the products.

\section{Main Results}
We introduce a non-anticipatory algorithm, dubbed iMWU, that achieves near optimality, while ensuring resource constraints are met in each time step. In more details, we develop a multi-phase version of the multiplicative weight update (MWU) algorithm, displayed in Algorithm \ref{alg:unknown_IB}. The iMWU algorithm crucially involves a virtual MWU process (displayed in Algorithm \ref{alg:IB}), which progressively produces weight vectors that trade-off among the amounts of rewards earned, resources consumed as well as the usage durations. The performance guarantee of iMWU crucially depends on the quantity $\xi = a_\text{max}/ c_\text{min}$, which can be interpreted as the maximum fraction of a resource allocated by an allocation decision in one time step. 
\begin{assumption}\label{ass:cap}
It holds that \textcolor{black}{$\xi \log(|\III_c| / \xi) \leq 1$}.
\end{assumption}
The assumption ensures that the DM has sufficient amounts of resources to buffer against the allocation error due to model uncertainty, and the stochastic variations due to the random outcomes. We remark that Assumption \ref{ass:cap} is strictly weaker than \citep{chen2017revenue,lei2020real},   who consider a scaling regime where $\frac{1}{\xi}$ grows with $T$. Instead, given that $|\III_c|$ is fixed, we only require $\frac{1}{\xi}$ to be larger than an absolute constant. The absolute constant \textcolor{black}{1} on the right hand side in Assumption \ref{ass:cap} is an artifact of our analysis, and it can potentially be replaced by a smaller constant by streamlining the constants involved in our analysis. In addition, we emphasize that Assumption \ref{ass:cap} is only needed for our theoretical analysis, and our proposed iMWU algorithm can still be implemented even when the Assumption is violated. In passing, we remark that Assumption \ref{ass:cap} is stronger than \citep{rusmevichientong2020dynamic,feng2020near}, who only require $\xi \leq 1$ in single reward maximization settings. These works consider the case of full model certainty, which allows them to carefully allocate their reusable resources based on the information about the customer arrival process. In our case, we require a higher capacity to buffer against errors due to model uncertainty.

The main result of this manuscript is the performance guarantee of our iMWU algorithm:

\begin{theorem}\label{thm:main}
Suppose Assumption \ref{ass:cap} holds. For any $\delta\in (0, 1)$, iMWU satisfies
\begin{equation}\label{eq:main}
\min_{i\in \III_r}\left\{\frac{1}{T}\sum^{T}_{t = 1} W_i(t)\right\} \geq \left(1 -3\sqrt{\xi \log \frac{|\III_c|}{\xi}}\right)\cdot \text{opt}(\text{IP-C}) - O\left(\frac{\log T}{\sqrt{T}}\right) 
\end{equation}
with probability at least \textcolor{black}{$1 - (8\log T)\delta$}. The $O(\cdot)$ notation hides a multiplicative factor in terms of $d_\text{max}, v_\text{max}, w_\text{max}$ and $\log((|\III_r| + |\III_c|) / \delta)$, but the factor is independent of $|\JJJ|, |\KKK|, T$.
\end{theorem}
The Theorem shows that the objective value of iMWU is closer to the optimum when $T$ and $c_{\text{min}}$ increase. A larger $c_{\text{min}}$, which leads to a smaller $\xi$, means that the DM is endowed with more resource units to buffer against the model uncertainty and random variations in the allocation outcomes. Remarkably, the performance bound does not degrade when the customer type set $\JJJ$ or the action set $\KKK$ grows. The independence to $|\JJJ|, |\KKK|$ are essential to assortment planning applications, where the size of $\JJJ$ is often exponential in dimension of the customers' feature vectors, and the size of $\KKK$ is often exponential in the maximum size of an assortment. While an optimal policy crucially depend on $\{p_j\}_{j\in \JJJ}$, we demonstrate that, in order to achieve near optimality, it is not necessary to estimate each and every $p_j$ accurately. Rather, our iMWU algorithm strives to achieve a near-optimal trade-off among the rewards earned, the resource units occupied and their usage durations, which only requires us to have an accurate estimate on certain weighted sum that involves $\{p_j\}_{j\in \JJJ}$. 

The performance guarantee in Theorem \ref{thm:main} consists of two components. The first component is the multiplicative factor $1 -3\sqrt{\xi \log(|\III_c|/\xi)}$, which is an approximation ratio. The approximation ratio captures the error due to the stochastic deviations of the random usage and rewards as compared to their expectations, as well as model uncertainty. As highlighted in the literature review, similar approximation ratios are derived in related works on online non-reusable resource allocations and offline reusable resource allocations with stochastic customer arrivals. The second component is the additive error term $O(\log T / \sqrt{T})$. The error term represents the estimation error, primarily due to the model uncertainty on $\textbf{p}, T$, as well as the fact that the collection $\{(w_{j(t), k}, v_{j(t), k})\}_{k\in \KKK}$ for customer $t$ is only revealed at his/her arrival. In the regime of large $T$, the error from the first component dominates the error from the second, which tends to 0 when $T$ grows.

Finally, we remark on the our result in relation to existing works on offline reusable resource allocations. In the case of single reward and single resource ($|\III_r| = |\III_c|=1$), our result leads to a $(1 -3\sqrt{\xi \log(|\III_c|/\xi)})$-approximation when $T$ is large, which nearly matches the $1 - O(\sqrt{\xi})$ approximation factor by \citep{levi2010provably} (recall that their approximation factor approaches $1 - \sqrt{2\xi/\pi}$ when $\xi$ tends to 0), in the sense that our loss term involves an additional log factor. We also remark that \citep{feng2020near,lei2020real} achieve approximation factors of $1 - O(\sqrt{\xi})$ in their respective settings. 
Nevertheless, it is worth emphasizing that our result does not assume the knowledge of $\{p_j\}_{j\in \JJJ}$, which corresponds to the arrival rates of each customer class in the settings in \cite{levi2010provably,feng2020near,lei2020real}. In a nutshell, our proposed policy achieves a nearly matching performance guarantee compared to \citep{levi2010provably,feng2020near,lei2020real} without assuming full model certainty. 

\subsection{An iterated MWU (iMWU) Algorithm}

Our main algorithm, the iterative MWU (iMWU) Algorithm, is presented in Algorithm \ref{alg:unknown_IB}. We partition the planning horizon into multiple phases $q=-1, 0, 1, \ldots$, where each phase $q$ consists of $\tau^{(q)}$ time steps. We set $\tau^{(-1)}=d_\text{max}$, and let $\tau^{(q)}$ doubles $\tau^{(q-1)}$ for all $q \geq 0$. For the default input values, we let 
\begin{equation}\label{eq:param}
\eta = \sqrt{\xi \log \frac{|\III_c|}{\xi}} ,\quad \bar{\epsilon}^{(q)}_D = \frac{\epsilon^{(q)}_D}{c_\text{min}}\text{, where }\epsilon^{(q)}_D = 8\gamma\sqrt{\frac{\log ((|\III_r|+|\III_c|) / \delta)}{\tau^{(q-1)}}}.
\end{equation}
Assumption \ref{ass:cap} ensures that $\eta \leq 1$. The discounting parameters $\eta$ and $\bar{\epsilon}^{(q)}_D$ serve to slim the customer flow, so that with high probability, the DM has enough available resource units to allocate over the planning horizon. Phase $-1$ \textcolor{black}{serves as a warm-up phase for the DM to learn about the latent model}.

In addition, the iMWU algorithm assumes the access to an
\emph{optimization oracle} $\kappa: \Delta^{\III_r\cup \III_c} \times \JJJ\rightarrow\KKK$, which has the property that for any $({\bm \phi}, {\bm \psi} ) \in \Delta^{\III_r\cup \III_c}$ and $j\in \JJJ$, the oracle returns
\begin{equation}\label{eq:opt_oracle}
\kappa({\bm \phi}, {\bm \psi} , j) \in \text{argmax}_{k \in \KKK} \left\{\sum_{i\in \III_r}\phi_{i}w_{ijk} - \sum_{i\in \III_c} \psi_{i} v_{ijk}\right\}.
\end{equation} 
In the case when $\KKK$ is small, for example when $\KKK$ only consists of the actions of ``accept'' or ``reject'' in the context of admission control, the optimization oracle $\kappa$ can be realized by evaluating the objective in (\ref{eq:opt_oracle}) for each action $k$. By contrast, in the applications involving assortment planning, the size of $\KKK$ could be large, since $\KKK$ consists of assortments, which are represented as subsets of products. In such applications, the size of $\KKK$ could be exponential in the number of products, which makes a direct enumeration infeasible. Nevertheless, the resulting assortment optimization problem in (\ref{eq:opt_oracle}) has been shown to be solvable in time polynomial on the number of products on a variety of choice models, a prominent example being the Multi-Nomial Logit choice model \citep{DavisGT13} under a cardinality constraint on the assortment size. Other examples include nested logit choice model under cardinality constraint \cite{DavisGT14}, and Markov chain choice model \cite{BlanchetGG16}. 

\begin{algorithm}[t]
\caption{Iterated MWU Algorithm}\label{alg:unknown_IB}
\begin{algorithmic}[1]
\State Input: Confidence parameter $\delta\in (0, 1)$.
\State Set error parameters $\{\bar{\epsilon}^{(q)}_D\}^\infty_{q=0}\subset (0, \infty)$ and discount parameter $\eta> 0$ according to (\ref{eq:param}).
\State Set $\tau^{(-1)}=d_\text{max}$, and set $\tau^{(q)}= 2\tau^{(q-1)}$ for all $q=0,1, \ldots$.
\For{time steps $t= 1, \ldots, \tau^{(-1)}$} \Comment{Phase $-1$.}
\State If each resource type $i\in \III_c$ has $\geq a_\text{max}$ units available, select an arbitrary action. Otherwise, select the null action $k_\textsf{null}$. 
\EndFor
\For{phases $q = 0, 1, \ldots$} 
\State Solve $\text{(LP-RS)}^{(q)}$ (see (\ref{eq:LP-RS-q})) for its optimal value $\hat{\lambda}^{(q)}_{*}$.
\State Invoke Algorithm \ref{alg:IB}, which returns the weight vector set $\Theta^{(q)} = \{ ({\bm \phi}^{(q)}(s),  {\bm \psi}^{(q)}(s)) \}^{\tau^{(q-1)}/2}_{s=1}$. \label{alg:Theta}
\For{time steps $t = \tau^{(q-1)}+1, \ldots, \tau^{(q)}$}
\State Sample a weight vector $(\boldsymbol{\phi}(t),\boldsymbol{\psi}(t))$ uniformly at random from $\Theta^{(q)}$.
\State Observe $j(t)$, and compute \label{alg:tilde}
\begin{equation*}
\tilde{k}(t)=\kappa(\boldsymbol{\phi}(t),\boldsymbol{\psi}(t), j(t)).
\end{equation*}
\State Set $\hat{k}(t) = \tilde{k}(t)$ with probability $\frac{1}{1 + \bar{\epsilon}^{(q)}_D + \eta}$ and $\hat{k}(t) = k_\textsf{null}$ with probability $1 - \frac{1}{1 + \bar{\epsilon}^{(q)}_D + \eta}$.\label{alg:hat}
\If{$\sum^{t-1}_{\tau=1} A_i(\tau)\mathbf{1}(D_i(\tau) \geq t - \tau+1)\leq c_i - a_\text{max}~\forall i\in \III_c$}
\State Take action $k(t) = \hat{k}(t)$.
\Else
\State Take action $k(t) = k_\textsf{null}$.
\EndIf
\EndFor
\EndFor
\end{algorithmic}
\end{algorithm}

\begin{algorithm}[h]
\caption{Virtual MWU process (invoked at the start of phase $q$)}\label{alg:IB}
\begin{algorithmic}[1]
\State Input: Observations during time steps $1, \ldots, \tau^{(q-1)}/2$, estimate $\hat{\lambda}^{(q)}_*$.
\State Initialize, for each $i\in \III_r, $ parameters $\phi^{(q)}_{i}(1) = \frac{1}{|\III_r| + |\III_c|}$ and $\Gamma^{(q)}_{i}(1) = 0$.
\State Initialize, for each $i\in \III_c, $ parameters $\psi^{(q)}_{i}(1) = \frac{1}{|\III_r| + |\III_c|}$ and $\Xi^{(q)}_{i}(1) = 0$.
\State Initialize $\eta(s) = \frac{\sqrt{\log(|\III_r|+|\III_c|)}}{\gamma \sqrt{s}}$ for each $s = 1, 2, \ldots, \tau^{(q-1)}/2$, where $\gamma = \max\{w_{\text{max}}, v_{\text{max}}\}$.
\For{virtual time steps $s = 1, \ldots, \tau^{(q-1)}/2$}
\State \label{alg:MWU_opt_oracle}Compute the \emph{virtual action}
\begin{equation*}
k^{(q)}(s)= \kappa\left({\bm \phi}^{(q)}(s), {\bm \psi}^{(q)}(s), j(s)\right)
\end{equation*}
by invoking the optimization oracle $\kappa$ in (\ref{eq:opt_oracle}).
\State Compute the exponent parameters
\begin{subequations}
\begin{alignat}{2}
&\Gamma^{(q)}_i(s+1)=\Gamma^{(q)}_i(s) + \left[w_{i, j(s), k^{(q)}(s)} -\left( \hat{\lambda}^{(q)}_{*} - \epsilon^{(q)}_C\right)\right] \quad &\forall i \in \III_r, \nonumber \\
&\Xi^{(q)}_{i}(s+1)=\Xi^{(q)}_{i}(s) + \left[ - v_{i, j(s), k^{(q)}(s)}+ \min\{c_i, v_\text{max}\}\right] \quad &\forall i \in \III_c, \nonumber
\end{alignat}
\end{subequations}
where $\epsilon^{(q)}_{C}$ is defined in (\ref{eq:epsilon_C}).
\State Compute the weights ${\bm \phi}^{(q)}(s+1) = \{\phi^{(q)}_i(s+1)\}_{i\in \III_r}$, ${\bm \psi}^{(q)}(s+1) = \{\psi^{(q)}_i(s+1)\}_{i\in \III_c}$:
\begin{subequations}
\begin{alignat}{2}
&\phi^{(q)}_i(s+1)=\frac{\exp\left[-\eta(s+1)\cdot \Gamma^{(q)}_i(s+1)\right]}{\sum_{i' \in \III_r} \exp\left[-\eta(s+1) \cdot \Gamma^{(q)}_{i'}(s+1)\right] +\sum_{i' \in \III_c} \exp\left[-\eta(s+1)\cdot \Xi_{i'}^{(q)}(s+1)\right]} \quad \forall i \in \III_r, \nonumber \\
&\psi_i^{(q)}(s+1) =\frac{\exp\left[-\eta(s+1) \cdot \Xi_i^{(q)}(s+1)\right]}{\sum_{i' \in \III_r} \exp\left[-\eta(s+1) \cdot \Gamma_{i'}^{(q)}(s+1)\right]+\sum_{i' \in \III_c} \exp\left[-\eta(s+1) \cdot \Xi_{i'}^{(q)}(s+1)\right]} \quad \forall i \in \III_c. \nonumber
\end{alignat}
\end{subequations}
\EndFor
\State Return the weight vector set $\Theta^{(q)} = \{ ({\bm \phi}^{(q)}(s),  {\bm \psi}^{(q)}(s)) \}^{\tau^{(q-1)}/2}_{s=1}$.
\end{algorithmic}
\end{algorithm}

Starting form phase 0, the iMWU Algorithm runs the following three steps for each phase:
\begin{itemize}
\item \emph{Estimate the value of $\lambda^E_*$} (line 8 of Algorithm \ref{alg:unknown_IB}): We use the observation collected during time steps $(\tau^{(q-1)} / 2) + 1, \ldots, \tau^{(q-1)}$, namely the second half of phase $q-1$, to compute \textcolor{black}{an estimate $\hat{\lambda}^{(q)}_*$ to $\lambda^E_*$.} The computation of $\hat{\lambda}^{(q)}_*$ first involves defining, for each $j\in \JJJ$, the empirical probability distribution
$$
\hat{p}^{(q)}_{j} = \frac{2}{\tau^{(q-1)}} \sum^{\tau^{(q-1)}}_{t = (\tau^{(q-1)} / 2) + 1} \mathbf{1}(j(t) = j).
$$
Recall that $v_{ijk} = \mathbb{E}[A_{ijk}D_{ijk}]$. The observation collected during time steps $(\tau^{(q-1)} / 2) + 1, \ldots, \tau^{(q-1)}$ is used to construct a sample average approximation (SAA) problem:
\begin{align}
\text{(LP-RS)}^{(q)}:  ~\max\limits_{\hat{y}_{jk}^{(q)}}  &~\hat{\lambda}^{(q)} & \label{eq:LP-RS-q}\\
\text{s.t.}  &\sum_{j \in \JJJ} \sum_{k\in \KKK} \hat{p}^{(q)}_{j} w_{ijk} \hat{y}_{jk}^{(q)} \geq \hat{\lambda}^{(q)}    &\quad &\forall i\in \III_r   \nonumber\\
&\sum_{j\in \JJJ}\sum_{k\in \KKK} \hat{p}^{(q)}_{j} v_{ijk} \hat{y}_{jk}^{(q)} \leq c_i       &\quad & \forall i\in \III_c \nonumber\\
&\sum_{k\in \KKK} \hat{y}_{jk}^{(q)} \leq 1      &\quad &\forall j\in \JJJ \nonumber\\
&\hat{y}_{jk}^{(q)} \geq 0      &\quad &\forall j\in \JJJ, ~k\in \KKK \nonumber.
\end{align}
to (LP-S), and we denote $\hat{\lambda}^{(q)}_*$ as the optimal value of $\text{(LP-RS)}^{(q)}$. It is worth noticing that $\text{(LP-RS)}^{(q)}$ is not a direct sample average approximation of (LP-E), where we only replace $p_j$ by its empirical value. In addition, we replace $\sum^t_{\tau=1} \mathbb{E}[A_{ijk} \mathbf{1}(D_{ijk} \geq t - \tau +1)]$ in (LP-E) by $v_{ijk}$ for any $t \in [T]$ in $\text{(LP-RS)}^{(q)}$. In this way, our SAA problem $\text{(LP-RS)}^{(q)}$ is more tractable than (LP-E). Indeed, the numbers of both the constraints and the variables of the former do not grow with $T$. In the assortment planning application when $\KKK$ could be of size exponential in $|\III_c|$, existing research works \cite{bront2009column} can be applied to solve $\text{(LP-RS)}^{(q)}$ efficiently via column generation. Moreover, the knapsack structure of $\text{(LP-RS)}^{(q)}$ allows us to construct an online algorithm that achieves a time average objective converging to $\hat{\lambda}^{(q)}_*$. The construction paves the way for the following virtual MWU process.
\item \emph{Run a virtual MWU process} (line 9 of Algorithm \ref{alg:unknown_IB}, Algorithm \ref{alg:IB}): The function of this step is to generate a set of weight vectors $\Theta^{(q)}$. These weight vectors are the dual variables generated by an MWU algorithm when solving the online feasibility problem of $\text{(LP-RS)}^{(q)}$ given $\hat{\lambda}^{(q)}_*$. Specifically, the MWU algorithm selects $k^{(q)}(s) \in \KKK$ in each time step $s = 1, \ldots, \tau^{(q-1)}/2$ assuming the optimization oracle $\kappa$ defined in (\ref{eq:opt_oracle}). 

\item \emph{Run a sampling process} (line 10 to 19 of Algorithm \ref{alg:unknown_IB}): The weight vectors generated in the previous step can be seen as the \textcolor{black}{penalty} of violating constraints of $\text{(LP-RS)}^{(q)}$ given different customer samples. Therefore, these weight vectors altogether describe the dual prices that trade-off between the rewards earned and the resources consumed. We uniformly sample a weight vector from the set, and take action $k(t)$ according to the sampled weight. This step also facilitates our analysis on bridging the gap between the resource constraints in $\text{(LP-RS)}^{(q)}$ and the actual resource consumption.
\end{itemize}


\subsection{Analysis}
In this section, we establish our performance guarantee on iMWU, by providing an overview on the proof of Theorem \ref{thm:main}.

\textbf{Linear Programming Benchmarks.}
We begin our analysis by providing a more tractable offline benchmark than the benchmark (LP-E). Firstly, recall that by Lemma \ref{lem:benchmark1}, we know that $\lambda^E_* \geq \mathbb{E}[\text{opt(IP-C)}]$, and therefore (LP-E) can serve as an offline benchmark for any online algorithm. However, 
it is hard to directly show the convergence of the time average reward to $\lambda^E_*$, while satisfying the resource constraints in each time step in an online setting. Indeed, the value of $T$, which is not known, crucially impact an optimal solution to (LP-E). To compensate for these drawbacks, in Algorithm \ref{alg:unknown_IB} we propose $\text{(LP-RS)}^{(q)}$, which is the sample average approximation of the following ``steady-state'' LP:
\begin{subequations}
\begin{alignat}{2}
\text{ (LP-S): }
\max\limits_{x_{jk}}  &~\lambda^S & \nonumber\\
\text{s.t.}  &\sum_{j \in \JJJ} \sum_{k\in \KKK} p_j w_{ijk} y_{jk}\geq \lambda^{S}     &\quad &\forall i\in \III_r   \nonumber\\
&\sum_{j\in \JJJ}\sum_{k\in \KKK} p_j v_{ijk} y_{jk} \leq c_i       &\quad & \forall i\in \III_c \nonumber\\
&\sum_{k\in \KKK} y_{jk} \leq 1      &\quad &\forall j\in \JJJ \nonumber\\
&y_{jk}\geq 0      &\quad &\forall j\in \JJJ, ~k\in \KKK \nonumber.
\end{alignat}
\end{subequations}
We denote the optimal solution of the steady-state LP as $y^*_{jk}$. To facilitate our forthcoming discussions, we abbreviate opt(LP-S), the optimal objective value of (LP-S), as $\lambda_*$. 
\begin{lemma}\label{lem:benchmark2}
$T \cdot\text{opt}(\text{LP-E})- d_{\max} w_{\max} \leq T \lambda_* \leq T \cdot \text{opt}(\text{LP-E}),$ where $\lambda_* = \text{opt}(\text{LP-S})$.
\end{lemma}
Lemma \ref{lem:benchmark2} is proved in Appendix Section \ref{pf:lemma_benchmark2}. The closeness of $\lambda_*$ to $\text{opt}(\text{LP-E})$ facilitates our analysis on the virtual MWU process. We further demonstrate that $T\cdot \text{opt(LP-E)} - T \lambda_* = \Omega(d_{\max})$ by a deterministic instance. Let $d$ be some positive even integer. We have $1$ reward, $1$ resource (with capacity $\frac{d}{2}$), $1$ customer type (with arrival rate $1$) and $2$ non-null actions ($\KKK=\{k_0,k_1,k_2\}$). With slight abuse of notation, for action $k_1$ (resp. action $K_2$), we let the reward earned, the resource consumed and the usage duration be $W_{k_1}=\frac{3}{4}$, $A_{k_1}=1$ and $D_{k_1}=\frac{d}{2}$ (resp. $W_{k_1}=1$, $A_{k_1}=1$ and $D_{k_1}=d$) with certainty. When $T=\frac{d}{2}$, the optimal solution of (LP-S) is to always take action $k_1$ which results in $T \lambda^S_*=\frac{d}{2} \frac{3}{4} = \frac{3d}{8}$. However, the optimal solution of (LP-E) is to take action $k_2$ in each time step and $T \lambda^E_*=\frac{d}{2} \cdot 1 = \frac{d}{2}$. In this case, $T \lambda^E_* - T \lambda^S_* = \frac{d}{8}$.

\textbf{Coupling Argument.} After replacing (LP-E) with a more tractable benchmark (LP-S), we focus on analyzing the performance of iMWU, by bounding the total type $i$ reward collected in phase $q$, for every $i\in \III_r$ and $q\in \{0, 1, \ldots, \lceil \log_2 (T/l)\rceil\}$. The analysis of the rewards collected within a phase is not straightforward. Indeed, the availability of each resource is intricately dependent on the amount of resources allocated in the previous time steps, as well as their usage durations. We disentangle the intricate dependence among the rewards $W_i(\tau^{(q-1)}+1), \ldots, W_i(\tau^{(q)})$ via a coupling argument. 

The argument goes as follows. Firstly, for each time step $t = \tau^{(q-1)}+1, \ldots, \tau^{(q)}$ in phase $q$, we denote $(\tilde{W}(t), \tilde{A}(t), \tilde{D}(t)) \sim {\cal O}_{j(t), \tilde{k}(t)}$, where the action $\tilde{k}(t)$ is computed in Line \ref{alg:tilde} in Algorithm \ref{alg:unknown_IB}. Next, we set $B(t)\sim \text{Bern}(1 / (1+\bar{\epsilon}^{(q)}_D + \eta))$, where $\text{Bern}(a)$ is the Bernoulli distribution with mean $a$. The random variables $\{B(t)\}^T_{t=1}$ are jointly independent, and they are independent of $\{(\tilde{W}(t), \tilde{A}(t), \tilde{D}(t))\}^T_{t=1}$. Then, the stochastic outcomes $(\hat{W}(t), \hat{A}(t), \hat{D}(t))\sim {\cal O}_{j(t), \hat{k}(t)}$ under the randomly chosen action $\hat{k}(t)$ (see Line \ref{alg:hat} in Algorithm \ref{alg:unknown_IB}) satisfy the following equalities with certainty:
\begin{align}
\hat{W}_i(t) &= B(t)\cdot \tilde{W}_i(t) \qquad\qquad\qquad\qquad \qquad\quad \text{ for all $i\in \III_r$},\nonumber\\
\hat{A}_i(t) &= B(t)\cdot \tilde{A}_i(t), \qquad \hat{D}_i(t) = B(t)\cdot \tilde{D}_i(t)~ \text{ for all $i\in \III_c$}.\label{eq:hatWAD}
\end{align}
By the definition of $k(t)$, we see that the actual outcomes $W_i(t), A_i(t), D_i(t)$ satisfy the following with certainty:
\begin{align}
W_i(t) &= \hat{W}_i(t)\cdot \mathbf{1}\left(\sum^{t-1}_{\tau=1} A_i(\tau)\mathbf{1}(D_i(\tau) \geq t - \tau+1)\leq c_i - a_\text{max}~\forall i\in \III_c\right) \text{ for all $i\in \III_r$},\nonumber\\
A_i(t) &= \hat{A}_i(t)\cdot \mathbf{1}\left(\sum^{t-1}_{\tau=1} A_i(\tau)\mathbf{1}(D_i(\tau) \geq t - \tau+1)\leq c_i - a_\text{max}~\forall i\in \III_c\right)\text{ for all $i\in \III_c$}, \label{eq:WAD}\\
D_i(t) &= \hat{D}_i(t)\cdot \mathbf{1}\left(\sum^{t-1}_{\tau=1} A_i(\tau)\mathbf{1}(D_i(\tau) \geq t - \tau+1)\leq c_i - a_\text{max}~\forall i\in \III_c\right) \text{ for all $i\in \III_c$}.\nonumber
\end{align}
The set of equations (\ref{eq:WAD}) fleshes out the (rather complex) dependence between the outcomes $W(t), A(t), D(t)$ in a time step $t$ and the outcomes in the previous time steps. The dependence is reflected by the indicator random variable on whether the DM is able to make an allocation at time $t$,  while ensuring that all resource constraints hold with certainty. 

We analyze the total type $i$ reward in a phase $q$ via the coupling. 
In the following, all equalities and inequalities hold almost surely.
\begin{align}
&\sum^{\tau^{(q)}}_{t = \tau^{(q-1)}+1} W_i(t) \nonumber\\
 = &\sum^{\tau^{(q)}}_{t = \tau^{(q-1)}+1} \hat{W}_i(t)\cdot \mathbf{1}\left(\sum^{t-1}_{\tau=1} A_i(\tau)\mathbf{1}(D_i(\tau) \geq t - \tau+1)\leq c_i - a_\text{max}~\forall i\in \III_c\right)\nonumber\\
\geq &\sum^{\tau^{(q)}}_{t = \tau^{(q-1)}+1} \hat{W}_i(t)\cdot \mathbf{1}\left(\sum^{t-1}_{\tau=1} \hat{A}_i(\tau)\mathbf{1}(\hat{D}_i(\tau) \geq t - \tau+1)\leq c_i - a_\text{max}~\forall i\in \III_c\right)\label{eq:couple1}\\
= &\sum^{\tau^{(q)}}_{t = \tau^{(q-1)}+1} \hat{W}_i(t)\cdot \mathbf{1}\left(\sum^{t-1}_{\tau=\max\{t-d_\text{max}, 1\}} \hat{A}_i(\tau)\mathbf{1}(\hat{D}_i(\tau) \geq t - \tau+1)\leq c_i - a_\text{max}~\forall i\in \III_c\right)\label{eq:couple2}\\
\geq &\sum^{\tau^{(q)}}_{t = \min\{\tau^{(q-1)}+1+d_\text{max}, \tau^{(q)}\}} \hat{W}_i(t)\cdot \mathbf{1}\left(\sum^{t-1}_{\tau=\max\{t-d_\text{max}, 1\}} \hat{A}_i(\tau)\mathbf{1}(\hat{D}_i(\tau) \geq t - \tau+1)\leq c_i - a_\text{max}~\forall i\in \III_c\right)\label{eq:couple3}.
\end{align}
Step (\ref{eq:couple1}) is because the coupling (\ref{eq:WAD}) ensures $A_i(t)\leq \hat{A}_i(t)$ and $D_i(t)\leq \hat{D}_i(t)$ almost surely. Step (\ref{eq:couple2}) uses the model assumption that $\hat{D}_i(t)\in \{0, \ldots, d_\text{max}\}$ almost surely. Finally, the inequality in (\ref{eq:couple3}) facilitates our analysis, since the time indexes associated with the stochastic outcomes $\hat{W}, \hat{A}, \hat{D}$ are all contained in $\{\tau^{(q-1)}+1, \ldots, \tau(q)\}$, the time steps in phase $q$. Conditioned on the weight vector set $\Theta^{(q)}$ (Line \ref{alg:Theta} in Algorithm \ref{alg:IB}), which is $\sigma({\cal H}(\tau^{(q-1)}))$-measurable, the random outcomes $\{(\hat{W}(t), \hat{A}(t), \hat{D}(t))\}^{\tau^{(q)}}_{t = \tau^{(q-1)}+1}$ defined in (\ref{eq:hatWAD}) are independently and identically distributed. Consequently, the random outcomes $\{(\hat{W}(t), \hat{A}(t), \hat{D}(t))\}^{\tau^{(q)}}_{t = \tau^{(q-1)}+1}$ are much easier to analyse than the set of actual random outcomes $\{(W(t), A(t), D(t))\}^{\tau^{(q)}}_{t = \tau^{(q-1)}+1}$. 

Next, we bound the sum (\ref{eq:couple3}) from below, which leads to the proof of Theorem \ref{thm:main}. To facilitate our discussion, we denote the shorthand $\bar{\tau}^{(q-1)}+1 = \min\{\tau^{(q-1)}+1+d_\text{max}, \tau^{(q)}\}$. For each $t\in \{\bar{\tau}^{(q-1)}+1, \ldots, \tau^{(q)}  \}$, we denote the Bernoulli random variable $$L(t) = \mathbf{1}\left(\sum^{t-1}_{\tau=\max\{t-d_\text{max}, 1\}} \hat{A}_i(\tau)\mathbf{1}(\hat{D}_i(\tau) \geq t - \tau+1)\leq c_i - a_\text{max}~\forall i\in \III_c\right).$$ Observe that 
$$
(\ref{eq:couple3}) = \sum^{\tau^{(q)}}_{t = \bar{\tau}^{(q-1)}+1} \hat{W}_i(t) L(t).
$$
The analysis on (\ref{eq:couple3}) is based on the consideration of $\hat{\lambda}^{(q)}_*$, an empirical estimate to $\lambda_*$ using second half of the data collected in phase $q-1$, as well as the Virtual MWU process in Algorithm \ref{alg:IB}, which returns the weight vector set $\Theta^{(q)}$ using the first half of the data collected in phase $q-1$. 
 In the following lemma, we bound the estimation error of $\lambda_*$. To this end, recall that $\gamma = \max\{w_\text{max}, v_\text{max}\}$.

\begin{lemma}\label{lem:emp}
For any $\delta\in (0, 1)$ and any $q\in \{1, 2, \ldots\}$, we have
\textcolor{black}{
\begin{align}
& \lambda_* \geq \hat{\lambda}^{(q)}_* -  \epsilon^{(q)}_C &\text{w.p. } 1 - \delta, \nonumber \\
&\hat{\lambda}^{(q)}_* \leq \hat{\lambda}^{(q)}_{*}  + \epsilon^{(q)}_A + \epsilon^{(q)}_B & \text{w.p. } 1 - 2\delta, \nonumber
\end{align}}
where the error parameters $\epsilon^{(q)}_A, \epsilon^{(q)}_B,\epsilon^{(q)}_C\in O(1/\sqrt{\tau^{(q-1)}})$ are
\textcolor{black}{
\begin{align}
\epsilon^{(q)}_A &= 2\sqrt{\frac{2\gamma}{c_{\text{min}}\tau^{(q-1)}}\log \frac{|\III_c| }{\delta} } + \frac{4\gamma}{c_{\text{min}}\tau^{(q-1)}}\log \frac{|\III_c| }{\delta},\label{eq:epsilon_A}\\
\epsilon^{(q)}_B &= 2w_{\max} \sqrt{\frac{ \log(|\III_r|/\delta)}{\tau^{(q-1)}}}, \label{eq:epsilon_B}\\
\epsilon^{(q)}_C &= \min\left\{2 w_{\max}\left[ \sqrt{\frac{2 \log(1/\delta)}{\tau^{(q-1)}} }+ \frac{2 \log(1/\delta)}{\tau^{(q-1)}}\right], w_\text{max}\right\}  \label{eq:epsilon_C}.
\end{align}}
\end{lemma}
Lemma \ref{lem:emp} is proved in Appendix Section \ref{pf:lemma_emp}. Next, the judicious choices of exponential parameters in Algorithm \ref{alg:IB} leads the following performance guarantee concerning the weight vector set $\Theta^{(q)}$:
\begin{lemma}\label{lem:iMWU}
Let $\delta \in (0, 1)$. The following inequalities hold:\textcolor{black}{
\begin{align}
&\frac{2}{\tau^{(q-1)}} \sum^{\tau^{(q-1)}/2}_{s=1} \sum_{j\in \JJJ}p_j w_{i,j,\kappa({\bm \phi}^{(q)}(s), {\bm \psi}^{(q)}(s), j)} &\geq \lambda_* -\epsilon^{(q)}_A-\epsilon^{(q)}_B-\epsilon^{(q)}_C- \epsilon^{(q)}_D   & \quad \forall i\in \III_r \quad \text{w.p. } 1 - 6\delta, \nonumber\\
&\frac{2}{\tau^{(q-1)}}
\sum^{\tau^{(q-1)}/2}_{s=1} \sum_{j\in \JJJ}p_j v_{i,j,\kappa({\bm \phi}^{(q)}(s), {\bm \psi}^{(q)}(s), j)} &\leq c_i + \epsilon^{(q)}_D & \quad \forall i\in \III_c \quad \text{w.p. } 1 - 3\delta, \nonumber
\end{align} 
where $\epsilon^{(q)}_D = 8\gamma\sqrt{\frac{\log ((|\III_r|+|\III_c|) / \delta)}{\tau^{(q-1)}}} = O(1/\sqrt{\tau^{(q)}})$.} 
\end{lemma} 
Lemma \ref{lem:iMWU} is proved in Appendix Section \ref{pf:lemma_iMWU}. Lemma \ref{lem:iMWU} is instrumental for demonstrating the following two statements on Lemma $L(t), \hat{W}_i(t)$:
\begin{lemma}\label{lem:bounding_reject} 
Suppose \textcolor{black}{$\tau^{(q-1)} \geq \frac{384\gamma^2}{c_\text{min}^2}\log\frac{|\III_r| + |\III_c|}{\delta}$.} With probability at least \textcolor{black}{$1-3\delta$}, the inequality
$$\mathbb{E}[L(t) ~|~ {\cal H}(\tau^{(q-1)})] \geq 1 - 2\sqrt{\xi}$$
holds simultaneously for all $t\in \{\bar{\tau}^{(q-1)}+1, \ldots, \tau^{(q)}  \}$.
\end{lemma} 
\begin{claim}\label{claim:bound_hat_W}
With probability at least \textcolor{black}{$1-6\delta$}, the inequality 
\begin{equation}\label{eq:bound_hat_W}
\mathbb{E}[\hat{W}_i(t) ~|~ {\cal H}(\tau^{(q-1)})  ] \geq \left(1 -\textcolor{black}{\sqrt{\xi \log \frac{|\III_c|}{\xi}}} \right) \lambda^* - \epsilon^{(q)}_A - \epsilon^{(q)}_B-\epsilon^{(q)}_C- \textcolor{black}{(1+\lambda^*)} \epsilon^{(q)}_D  .
\end{equation}
holds simultaneously for all $t\in \{\bar{\tau}^{(q-1)}+1, \ldots, \tau^{(q)}  \}$ and all $i\in \III_r$. 
\end{claim}
To this end, we remark that all the error terms $\epsilon^{(q)}_A, \epsilon^{(q)}_B, \epsilon^{(q)}_C,  2\epsilon^{(q)}_D\in O(1/\sqrt{\tau^{(q)}})$. Consequently, as the horizon $T$ grows, the performance loss is dominated by the approximation factor \textcolor{black}{$\left(1 -\sqrt{\xi \log \frac{|\III_c|}{\xi}} \right)$}. Lemma \ref{lem:bounding_reject} is proved in Appendix Section \ref{pf:lema_bounding_reject}, and the proof uses Assumption \ref{ass:cap}. To illustrate the use of Lemma \ref{lem:iMWU}, we provide a proof of Claim \ref{claim:bound_hat_W} in the following:
\begin{proof}{Proof of Claim \ref{claim:bound_hat_W}}
By the coupling argument that constructs $\hat{W}_i(t)$, for every $t\in \{\bar{\tau}^{(q-1)}+1, \ldots, \tau^{(q)}  \}$ and every $i\in \III_r$, with certainty we have
\begin{align*}
\mathbb{E}[\hat{W}_i(t) ~|~ {\cal H}(\tau^{(q-1)})  ] & = \mathbb{E}[B(t)\tilde{W}_i(t) ~|~ {\cal H}(\tau^{(q-1)})  ] \nonumber\\
& = \mathbb{E}[B(t)]\cdot \mathbb{E}[\tilde{W}_i(t) ~|~ {\cal H}(\tau^{(q-1)})  ]\nonumber\\
& = \frac{1}{1+\bar{\epsilon}^{(q)}_D + \eta} \cdot \frac{2}{\tau^{(q-1)}} \sum^{\tau^{(q-1)}/2}_{s=1} \sum_{j\in \JJJ}p_j w_{i,j,\kappa({\bm \phi}^{(q)}(s), {\bm \psi}^{(q)}(s), j)} \nonumber.
\end{align*}
Applying Lemma \ref{lem:iMWU}, we know that the inequality 
\begin{align}
\frac{1}{1+\bar{\epsilon}^{(q)}_D + \eta} \cdot\frac{2}{\tau^{(q-1)}} \sum^{\tau^{(q-1)}/2}_{s=1} \sum_{j\in \JJJ}p_j w_{i,j,\kappa({\bm \phi}^{(q)}(s), {\bm \psi}^{(q)}(s), j)} \geq & \frac{1}{1+\bar{\epsilon}^{(q)}_D + \eta} \cdot \lambda_* -\epsilon^{(q)}_A-\epsilon^{(q)}_B-\epsilon^{(q)}_C- \epsilon^{(q)}_D \nonumber \\   
\geq & \lambda_* (1 - \eta) -\epsilon^{(q)}_A-\epsilon^{(q)}_B-\epsilon^{(q)}_C- (1 + \lambda_*)\epsilon^{(q)}_D \nonumber
\end{align}
holds for all $i\in \III_r$ with probability at least \textcolor{black}{$1-6\delta$}. Finally, unravelling the definition of $\eta$ leads to the desired lower bound. 
\hfill $\square$
\end{proof}
Combining Lemma \ref{lem:bounding_reject} and Claim \ref{claim:bound_hat_W}, we know that with probability at least $1-7\delta$, the inequality 
\begin{equation}\label{eq:bound_hat_W_and_L}
\mathbb{E}[\hat{W}_i(t) L(t) ~|~ {\cal H}(\tau^{(q-1)})  ] \geq \left(1 - \textcolor{black}{3\sqrt{\xi \log \frac{|\III_c|}{\xi}}} \right) \lambda^* - \epsilon^{(q)}_A - \epsilon^{(q)}_B-\epsilon^{(q)}_C- \textcolor{black}{(1+\lambda^*)} \epsilon^{(q)}_D  .
\end{equation}
holds simultaneously for all $t\in \{\bar{\tau}^{(q-1)}+1, \ldots, \tau^{(q)}  \}$ and all $i\in \III_r$. Finally, we use (\ref{eq:bound_hat_W_and_L}) to provide a lower bound to the total type $i$ rewards earned during phase $q$.
\begin{lemma}\label{lem:correlated}
With probability at least \textcolor{black}{$1-8\delta$}, the inequality $$\sum^{\tau^{(q)}}_{t=\tau^{(q-1)}+1}\hat{W}_i(t)L(t) \geq \left(1 -\textcolor{black}{3\sqrt{\xi \log \frac{|\III_c|}{\xi}}}\right) \lambda^* (\tau^{(q)} - \tau^{(q-1)}) - O\left(\sqrt{\tau^{(q)} -\tau^{(q-1)}}\right)$$
holds for all $i\in \III_r$. 
\end{lemma}
Lemma \ref{lem:correlated} is proved in Appendix Section \ref{pf:lemma_correlated}. By summing over $q\in \{-1, \ldots, \lceil \log_2 \frac{T}{\tau^{(-1)}}\rceil \}$ and a union bound over the phases, Theorem \ref{thm:main} is proved. We conclude our discussion by highlighting that the proof of Lemma \ref{lem:correlated} requires a conditional Chernoff inequality, stated below: 
\begin{lemma}[Conditional Multiplicative Chernoff Inequality]\label{lem:multi_AZ}
Suppose random variables $\{X_t\}^N_{t=1}$ satisfy the following properties:
\begin{enumerate}
\item $X_1, \ldots, X_N$ are jointly independent conditional on a $\sigma$-algebra ${\cal F}$,
\item $\Pr(X_t \in [0, B]\text{ for all $t\in \{1, \ldots, N\}$}) = 1$ for some $B\in \mathbb{R}_{>0}$.
\item There exists real numbers $\delta_{-}\in [0, 1]$ and $\mu_{-}\in \mathbb{R}_{> 0}$ such that
\begin{equation}\label{eq:mu_plus_minus_mini}
 \Pr\left(\mathbb{E}\left[\sum^N_{t=1} X_t\mid {\cal F}\right] < \mu_{-}\right) \leq \delta_-.
\end{equation}
\end{enumerate} 
Then the following concentration inequalities hold for any fixed but arbitrary $\delta\in (0,1)$:
\begin{align}
\Pr\left(\sum^N_{t=1} X_t < \mu_- - \sqrt{2B\mu_- \log\frac{1}{\delta}}\right) &\leq \delta + \delta_- .\label{eq:multi_AZ_minus_mini} 
\end{align}
\end{lemma}
A generalization of Lemma \ref{lem:multi_AZ} and its proof are provided in Appendix \ref{pf:lemma_multi_AZ}. The generalization includes both high probability upper and lower bounds to $\sum^{N}_{t=1}X_t$. While we only require the upper bound shown in the Lemma, the generalization could be of independent interest. We remark that Lemma \ref{lem:correlated} does not follow from a direct application of the conditional Chernoff inequality in Lemma \ref{lem:multi_AZ}, since the complicated correlation among $L(t)\hat{W}_i(t)$'s forbids a direct construction of martingale difference sequence. 

\section{Numerical experiments}\label{sec:num}
In our numerical experiments, we consider an assortment planning problem, with a set up in line with the description in \textbf{Assortment Planning} in Section \ref{sec:app}. One unit of resource $i$ is associated with a fixed price $r_i$. The DM influences the customers' demands through offering different assortments of resources. Contingent upon the arrival of a customer, say of type $j$, the DM decides the assortment $k \in \KKK$ to display, where $\KKK$ is a collection of subsets of $\III_c$. Let $q_{ijk}$ denote the probability for customer type $j$ to choose product $i$ when offered assortment $k$, where $i$ either belongs to $k$ or $i=0$, the no purchase option. If the customer chooses product $i$, one resource unit of $i\in \III_c$ is consumed, consequently we have $A_{ijk} \sim \text{Bern}(q_{ijk})$. The usage durations are defined in the same way as the discussions in  \textbf{Assortment Planning}, where $D_{ijk} = A_{ijk}\cdot \textsf{Duration}_{ij}$, and $A_{ijk},\textsf{Duration}_{ij}$ are independent. To facilitate our discussions, we denote $\textsf{d}_{ij} = \mathbb{E}[\textsf{Duration}_{ij}]$.


To define the reward types, we first partition the set of resources $\III_c$ is partitioned into 2 categories $\III^{(1)}_c$ and $\III^{(2)}_c$, meaning that $\III_c = \III^{(1)}_c \cup \III^{(2)}_c$ and $\III^{(1)}_c \cap \III^{(2)}_c = \emptyset$. The DM aims to simultaneously maximize the market share of each category, as well as maximize the total profit. For these different rewards, we introduce an array of KPI parameters $\boldsymbol{\sigma} = \{\sigma_i\}_{i \in \III_r}$ for normalization. Altogether, the multi-objective optimization problem is formulated as maximizing the minimum of three objectives indexed by $\III_r = \{1,2,3\}$: index 1 corresponds to the revenue: $W_{1,jk}=\frac{\sum_{i \in \III_c} r_i A_{ijk}}{\sigma_1}$, index 2 corresponds to the sales volume of category 1 resources: $W_{2,jk}=\frac{\sum_{i \in \III^{(1)}_c} A_{ijk}}{\sigma_2}$, and index 3 corresponds to the total sales volume of category 2 resources: $W_{3,jk}=\frac{\sum_{i \in \III^{(2)}_c} A_{ijk}}{\sigma_3}$. Note that these three objectives are normalized by $\boldsymbol{\sigma}$.

The choice probability $q_{ijk}$ is modeled by the multinomial logit (MNL) choice model as discussed in Section \ref{sec:app}. Each resource $i\in \III_c$ is associated with a feature vector $\boldsymbol{f}_i \in \mathbb{R}^m$, and each customer type $j\in \JJJ$ is associated with a set of feature vectors $\{\boldsymbol{b}_{ij} \}_{i \in \III_c}$, where $\boldsymbol{b}_{ij}  \in \mathbb{R}^m$ for each $i\in \III_c$. The feature vector $\boldsymbol{f}_i$ could involve the fixed price $r_i$. Under the MNL choice model, the choice probability is mathematically defined as
\begin{equation}\label{eq:mnl}
q_{ijk} = \frac{\exp(\boldsymbol{b}_{ij}^\top \boldsymbol{f}_i)}{1+\sum_{\ell \in k} \exp(\boldsymbol{b}_{\ell j}^\top \boldsymbol{f}_\ell)}
\end{equation}
if $i\in k$, and $q_{ijk} = 0$ if $i\not \in k$. In complement, the probability of no purchase is  $\frac{1}{1+\sum_{\ell \in k} \exp(\boldsymbol{b}_{\ell j}^\top \boldsymbol{f}_\ell)}$. 

To facilitate model transformation, we define ${\bm \sigma'} = \{\sigma'_i\}_{i \in \III_c}$ and ${\bm \phi'} = \{\phi'_i\}_{i \in \III_c}$ where
\begin{align}
\sigma'_i &= 
\begin{cases}
\sigma_2& \quad \text{for all }i \in \III^{(1)}_c, \\
\sigma_3& \quad \text{for all } i \in \III^{(2)}_c,
\end{cases}\nonumber\\
\phi'_i &= 
\begin{cases}
\phi_2& \quad \text{for all }i \in \III^{(1)}_c, \\
\phi_3& \quad \text{for all } i \in \III^{(2)}_c.
\end{cases}\nonumber
\end{align}
Consequently, the optimization oracle (see Equation (\ref{eq:opt_oracle})) can be expressed as
\begin{subequations}
\begin{alignat}{2}
\kappa({\bm \phi}, {\bm \psi} , j) \in & \text{argmax}_{k \in \KKK} \left\{\sum_{i\in \III_r}\phi_{i} w_{ijk} - \sum_{i\in \III_c} \psi_{i} v_{ijk}\right\}\nonumber\\
= &  \text{argmax}_{k \in \KKK} \left\{\left(\sum_{i \in \III_c} \frac{r_i \cdot \phi_1}{\sigma_1} + \sum_{i \in \III^{(1)}_c} \frac{\phi_2}{\sigma_2} +
\sum_{i \in \III^{(2)}_c} \frac{\phi_3}{\sigma_3} - \sum_{i \in \III_c} \textsf{d}_{ij} \cdot \psi_i \right) \cdot q_{ijk} \right\} \nonumber \\
= &  \text{argmax}_{k \in \KKK} \left\{\sum_{i \in \III_c} \left( \frac{r_i \cdot \phi_1}{\sigma_1} + \frac{\phi'_i}{\sigma'_i} - \textsf{d}_{ij} \cdot \psi_i\right) \cdot q_{ijk} \right\}. \nonumber
\end{alignat}
\end{subequations}
While the size of the action set $\KKK$ scales exponentially with the number of products, 
in the case of the MNL choice model (see (\ref{eq:mnl})), an optimal action $\kappa({\bm \phi}, {\bm \psi}, j)$ can be computed efficiently by solving the following LP \citep{davis2013assortment}:
\begin{align}
\max  &~\sum_{i \in \III_c} \left(\frac{r_i \cdot \phi_1}{\sigma_1} + \frac{\phi'_i}{\sigma'_i} - \textsf{d}_{ij} \cdot \psi_i \right) z_i & \nonumber\\
\text{s.t.}  & \sum_{i \in \III_c} z_i + z_0 = 1   \nonumber\\
& \sum_{i \in \III_c} \frac{z_i}{\exp(\boldsymbol{b}_{ij}^\top \boldsymbol{f}_i)} \leq n z_0 \nonumber\\
& 0 \leq \frac{z_i}{\exp(\boldsymbol{b}_{ij}^\top \boldsymbol{f}_i)} \leq z_0 \quad \forall i\in \III_c \nonumber
\end{align}
where the decision variables are $\{z_i: i \in \III_c \cup \{0\}\}$. 

We run simulations of iMWU and compare the results with the benchmark (LP-S), whose optimum value $\lambda_* $ satisfies  
\begin{equation}\label{eq:benchmark_compare}
\lambda_* \geq \text{opt(LP-E)}-\frac{d_\text{max}}{T} \geq \text{opt(IP-C)}-\frac{d_\text{max}}{T} 
\end{equation}
by Lemma \ref{lem:benchmark2}. We consider a synthetic dataset with $14$ types of resources indexed by $\III_c=\{1,2,\dots,14\}$ and $1000$ types of customers indexed by $\JJJ=\{1,2,\ldots,1000\}$. We allow offering any assortment of fewer than or equal to 5 products, and hence the action set $\KKK$ is of cardinality $\sum^5_{i=1} C^{14}_i=3472$. Each usage duration $\textsf{Duration}_{ij}$ follows a randomly generated discrete probability distribution, where the realized duration is upper bounded by $\frac{T}{5}$. To solve (LP-S) with a large action set comprising of assortments, we use the column generation technique introduced in \citet{bront2009column}.
We test two cases where $\xi=\frac{1}{20}$ and $\xi=\frac{1}{200}$ with $10$ simulations run for each. The results are illustrated in Figure 1 where each curve represent the average over 10 simulations, and the shaded area around each curve marks the variance over the simulations.

\begin{figure}[htbp]
\centering
\subfloat[Rewards losses $\xi=1/20$]{
\begin{minipage}[t]{0.5\textwidth}
\centering
\includegraphics[width=3in]{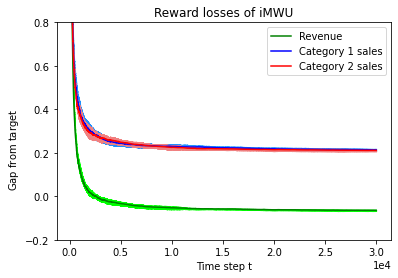}
\end{minipage}
}
\subfloat[Average normalized reward $\xi=1/20$]{
\begin{minipage}[t]{0.5\textwidth}
\centering
\includegraphics[width=3in]{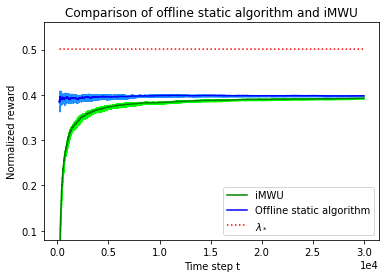}		    \end{minipage}
} \\
\subfloat[Rewards losses $\xi=1/200$]{
\begin{minipage}[t]{0.5\textwidth}
\centering
\includegraphics[width=3in]{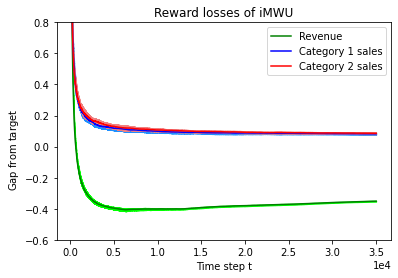}
\end{minipage}
}
\subfloat[Average normalized reward $\xi=1/200$]{
\begin{minipage}[t]{0.5\textwidth}
\centering
\includegraphics[width=3in]{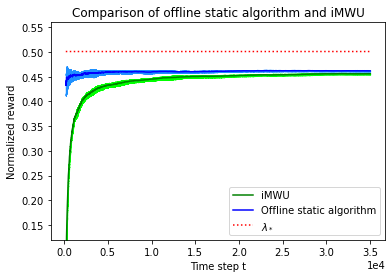}		    \end{minipage}
}
\caption{Performance comparison of different resource capacities}\label{fig:1}
\end{figure}
    
Figures \ref{fig:1}(a) and \ref{fig:1}(c) show how the 
\begin{equation}\label{eq:reward_gap}
\text{Reward Gap}_{i, t} = \frac{t  \lambda_* - \sum^t_{\tau=1} W_i(\tau)}{t  \lambda_*}
\end{equation}
for each objective $i \in \III_r = \{1, 2, 3\}$ varies as the time index $t$ grows under iMWU, for the cases of  $\xi=1/20$ and $\xi=1/200$ respectively.
The quantity $\lambda_*$ is the optimal value of the benchmark (LP-S). 
The rewards gaps can be negative. Indeed $\lambda_*$ is the minimum among the three per-time-step target rewards after normalization, and in both Figures \ref{fig:1}(a) and \ref{fig:1}(c) the normalized revenue is shown to be above $\lambda_*$, while the time average sales volumes of the two category are below $\lambda_*$. Consequently, the sales volumes are the binding reward constraints. When $\xi=1/20$, the sales losses for both categories converge to roughly $22\%$ of (LP-S)'s optimal value as $t$ gets larger; and when $\xi=1/200$, the sales losses converge to around $7\%$ for both categories. However, the revenue gaps are negative for both $\xi$ values since their values exceed $t \sigma_1 \lambda_*$ as $t$ gets larger. 

Figures \ref{fig:1}(b) and \ref{fig:1}(d) compare the per-time-step normalized reward of iMWU (see the forthcoming (\ref{eq:norm_iMWU}) with that of an offline static algorithm assuming known $\boldsymbol{p}$, derived from \citet{levi2010provably, devanur2019near}. The offline static algorithm (OSA) involves solving (LP-S) for an optimal solution $\{y^*_{jk}\}_{j\in \JJJ, k\in\KKK}$. The OSA selects action $k$ contingent upon a type-$j$ customer arrival with probability $y^*_{jk}/(1+\bar{\eta})$, where $\bar{\eta} = \Omega(\sqrt{\xi})$ is a discount parameter buffering against inventory shortages. It is worth mentioning that in the case of $|\III_c| = |\III_r| = 1$, \citet{levi2010provably} serves each type-$j$ customer with probability $y^*_{jk}$ as long as the available resources are sufficient for the allocation decision. Otherwise, the null action is chosen. They manage to achieve a $1 - O(\sqrt{\xi})$ approximation ratio when $\xi$ tends to 0. Nevertheless, since our model considers multiple resources and rewards, their analysis does not apply. We therefore involve a discount parameter in a similar manner as \cite{devanur2019near}. For evaluating algorithms, we define
\begin{align}
\text{Normalized Reward}^{OSA}(t) & = \min_{i \in \III_r} \left\{\frac{1}{t} \sum^t_{\tau=1} \sum_{k \in \KKK} W_{i,j(\tau),k}(\tau) \cdot X^{OSA}_k(\tau)\right\}, \nonumber \\
\text{Normalized Reward}^{iMWU}(t) & = \min_{i \in \III_r} \left\{\frac{1}{t} \sum^t_{\tau=1} \sum_{k \in \KKK} W_{i,j(\tau),k}(\tau) \cdot X^{iMWU}_k(\tau)\right\}. \label{eq:norm_iMWU}
\end{align}
The OSA achieves asymptotic performance guarantee of $(1 - \tilde{O}(\sqrt{\xi})) \lambda_*$. It is evident that iMWU, which suffers from model uncertainty, converges to the results of OSA, which has access to all model parameters, for both $\xi$ values as $t$ gets larger. Comparing the instances of $\xi=1/20$ and $\xi=1/200$, it is evident that iMWU (and OSA) perform better when the reciprocal of the relative resource capacity $\xi$ is smaller, which can be interpreted as the case when the DM is endowed with more resources. These results are consistent with our theoretical result on the multiplicative approximation factor $(1 - 3\sqrt{\xi \log (|\III_c|/ \xi)})$, which is closer to 1 when $\xi$ decreases. 

The difference between $\lambda^*$ and the normalized rewards under iMWU only serves as an \emph{upper bound} to the optimality gap $$\text{opt(IP-C)} - \text{Normalized Reward}^{iMWU}(t)$$ to within an additive error of $d_\text{max}/T$, by (\ref{eq:benchmark_compare}). It is not known if $\lambda^*$ or  $\text{Normalized Reward}^{iMWU}(t)$ is closer to the true optimum $\text{opt(IP-C)}$. The former case implies that the plots do illustrate the optimality gaps, while the latter case implies that $\lambda^*$ is a loose upper bound. In fact, as remarked by related works \citep{levi2010provably,feng2021online}, it is believed that solving opt(IP-C) is computationally hard even under full model uncertainty. 
While our theoretical result implies that the optimality gap tends to zero as $c_\text{min}, T$ grow, ascertaining a tighter bound on the optimality gap in the context of bounded $c_\text{min}$ requires a better approximation to $\text{opt(IP-C)}$.

\begin{figure}[H]
\centering
\includegraphics[width=0.5\textwidth]{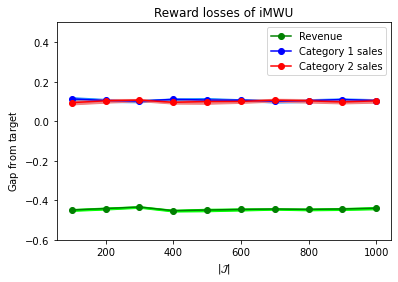}
\caption{Performance of different $\JJJ$ sizes}\label{fig:j}
\end{figure}

Finally, in Figure \ref{fig:j}, we further validate that our algorithmic performance is independent of the $\JJJ$ size. We fix $T=10000$, $\xi=1/200$ in Figure \ref{fig:j}. A family of instances is generated where $|\JJJ|$ varies between $100$ and $1000$. It can be seen that the Reward Gap (as defined in (\ref{eq:reward_gap})) of our iMWU algorithm are almost identical with different sizes of $\JJJ$.

\section{Conclusion}
In this paper, we develop an online algorithm for addressing reusable resource allocation problems. Our model aims to maximize multiple rewards generated by heterogeneous customers, whose arrival process model is not known. We propose an iterated multiplicative weight update (iMWU) algorithm achieving near-optimal performance guarantees. We remark that our model and algorithm capture a wide range of applications in Revenue Management, including admission control, pricing and assortment planning with reusable resources. Finally, we conduct numerical experiments to validate our theoretical results. For future research, it is interesting to consider the more challenging settings with non-stationary customer arrival processes. In addition, we believe that our the technical tools developed for iMWU could find applications in other settings involving reusable resources. 

\bibliographystyle{apalike}
\bibliography{reference}
\newpage
\begin{APPENDICES}
\section{Auxiliary Results}

\begin{proposition}[Azuma-Hoeffding Inequality]\label{prop:hoeffding}
Let $N$ be a positive integer and $B$ be a positive real number. Suppose the random variables $X_1, \ldots, X_N$ constitute a martingale difference sequence with respect to the filtration $\{{\cal F}_n\}^N_{n=0}$, i.e. $\mathbb{E}[X_n | {\cal F}_{n-1}] =0$ almost surely for every $n\in \{1, \ldots, N\}$. In addition, suppose $|X_n|\leq B$ almost surely for every $n\in \{1, \ldots, N\}$. For any $\delta \in (0, 1)$, it holds that
$$
\Pr\left[ \left| \frac{1}{N}\sum^N_{n=1} X_n \right| \leq B\sqrt{\frac{2\log (2/\delta)}{N}} \right] \geq 1 - \delta.
$$\hfill $\square$
\end{proposition}

\begin{proposition}[Multiplicative Chernoff inequality]\label{prop:multi_Chernoff}
Suppose random variables $\{X_t\}^N_{t=1}$ are independent, and that $\Pr(X_t \in [0, B]\text{ for all $t\in \{1, \ldots, N\}$}) = 1$ for some $B\in \mathbb{R}_{>0}$. Denote $\mu = \mathbb{E}[\sum^N_{t=1} X_t]$. The following concentration inequalities hold for any fixed but arbitrary $\delta\in (0,1)$:
\begin{align}
\Pr\left(\sum^N_{t=1} X_t > \mu + 2\sqrt{B\mu \log\frac{1}{\delta}} + 2B \log\frac{1}{\delta} \right) &\leq \delta, \label{eq:multi_chernoff_+}\\
\Pr\left(\sum^N_{t=1} X_t < \mu - \sqrt{2B\mu \log\frac{1}{\delta}}\right) &\leq \delta.\label{eq:multi_chernoff_-} 
\end{align}
\end{proposition}
While Proposition \ref{prop:multi_Chernoff} is well known, it can also be recovered by specializing Lemma \ref{lem:multi_AZ} via setting ${\cal F} = \emptyset$ and $\mu_+ = \mu_- = \mathbb{E}[\sum^N_{t=1} X_t] $ so that the inequalities in (\ref{eq:mu_plus_minus}) hold with $\delta_- = \delta_+ = 0$.

\begin{proposition}[Multiplicative Weight Update]\label{prop:MWU}
Let $\{\ell(s)\}^{\tau}_{s=1}$ be an arbitrary sequence of vectors, where $\ell(s) = (\ell_i(s))_{i\in \{\III_r, \cup\III_c\}}\in [-B, B]^{\III_r\cup \III_c}$ for each $s\in \{1, \ldots, \tau\}$. 
Consider the sequence of vectors $\vartheta(1), \ldots, \vartheta(\tau)$, where $\vartheta(s) = (\vartheta_{i}(s))_{i\in \III_r\cup \III_c}\in \Delta^{\III_r\cup\III_r}$ is defined as  
\begin{equation}\label{eq:vartheta}
\vartheta_{i}(s) = \frac{\exp\left[-\eta(s) \sum^{s-1}_{n=1}\ell_i(n) \right]}{\sum_{\iota\in \III_r\cup\III_c} \exp\left[-\eta(s) \sum^{s-1}_{n=1}\ell_\iota(n)\right] }, \text{ and }\eta(s) = \frac{\sqrt{\log(|\III_r| + |\III_c|)}}{B\sqrt{s}}
\end{equation}
 for each $s\in \{1, \ldots \tau\}, i\in \III_r\cup\III_c$. Then, for any $i\in \III_r\cup\III_c$, it holds that
$$
\frac{1}{\tau}\sum^\tau_{s=1} \ell_i(s) \geq \frac{1}{\tau}\sum^\tau_{s=1} \sum_{\iota\in \III_r\cup\III_c} \vartheta_{\iota}(s) \ell_{\iota}(s) - 2B\sqrt{\frac{\log(|\III_r| + |\III_c|)}{\tau} }. 
$$ \hfill $\square$
\end{proposition}
The proof of Proposition \ref{prop:MWU} can be found in Chapter 7.5 in \cite{orabona2019modern}.
\section{Proofs}
\subsection{Proof of Lemma \ref{lem:benchmark1}}\label{pf:lemma_benchmark1}
Let $\pi$ be a non-anticipatory feasible policy that achieves the expected optimum $\mathbb{E}[\lambda^C_*]$ in (IP-C), i.e. $
\mathbb{E}\left[\min_{i\in \III_r}\left\{\frac{1}{T}\sum^T_{t=1} \sum_{k\in \KKK} W_{i, j(t), k}(t) X^{\pi}_{k}(t)\right\}\right]= \mathbb{E}[\text{opt(IP-C)}].
$
Define $x = \{x_{jk}(t)\}_{j\in \JJJ, k\in \KKK, t\in \{1, \ldots, T\}}$ as $x_{jk}(t) = \Pr(X^\pi_k(t) = 1 | j(t) = j)$. We claim that $x$ is feasible to (LP-E), with objective value equal to $
\min_{i\in \III_r}\left\{\mathbb{E}\left[\frac{1}{T}\sum^T_{t=1} \sum_{k\in \KKK} W_{i, j(t), k}(t) X^{\pi}_{k}(t)\right]\right\}$, which is larger than or equal to $
\mathbb{E}\left[\min_{i\in \III_r}\left\{\frac{1}{T}\sum^T_{t=1} \sum_{k\in \KKK} W_{i, j(t), k}(t) X^{\pi}_{k}(t)\right\}\right]$. Thus, verifying the claims about the feasibility and the objective value proves the claim. 

We first verify the feasibility to (LP-E). Since the policy $\pi$ satisfies the reusable resource constraints, the inequality $\sum^t_{\tau = 1} \sum_{k\in \KKK} \mathbf{1}\left( D_{ij(\tau)k}\left(\tau\right) \geq t - \tau + 1\right) A_{ij(\tau)k}\left(\tau\right) X^{\pi}_{k}\left(\tau\right) \leq c_i$ holds for all $i\in \III_c, t\in \{1,\ldots, T\} $. Taking expectation (which is over over $X^{\pi}_k(\tau)$, $ D_{ijk}(\tau)$, $A_{ij(\tau)k}(\tau)$ and $j(\tau)$) for $\tau=1,\ldots,t$ on the left hand side gives
\begin{align}
& \mathbb{E}\left[\sum^t_{\tau = 1} \sum_{k\in \KKK} \mathbf{1}\left( D_{ij(\tau)k}\left(\tau\right) \geq t - \tau + 1\right) A_{ij\left(\tau\right)k}\left(\tau\right) X^{\pi}_{k}\left(\tau\right)\right] \nonumber \\
=& \sum^t_{\tau = 1} \sum_{k\in \KKK} \mathbb{E}\left[ \mathbb{E}\left[ \mathbf{1}\left( D_{ij(\tau)k}\left(\tau\right) \geq t - \tau + 1\right) A_{ij\left(\tau\right)k}\left(\tau\right) \mid j(\tau), X^{\pi}_{k}(\tau)\right] X^{\pi}_{k}(\tau) \right] \nonumber \\
=& \sum^t_{\tau = 1} \sum_{k\in \KKK} \mathbb{E}\left[ \mathbb{E}\left[\mathbb{E}\left[ \mathbf{1}\left( D_{ij(\tau)k}\left(\tau\right) \geq t - \tau + 1\right) A_{ij\left(\tau\right)k}\left(\tau\right) \mid j(\tau), X^{\pi}_{k}(\tau)\right] X^{\pi}_{k}(\tau)\mid j(\tau)\right] \right] \nonumber \\
=& \sum^t_{\tau = 1} \sum_{k\in \KKK} \sum_{j\in \JJJ} p_j  \mathbb{E}\left[ \mathbb{E}\left[ \mathbf{1}\left( D_{ij(\tau)k}\left(\tau\right) \geq t - \tau + 1\right) A_{ij\left(\tau\right)k}\left(\tau\right) \mid j(\tau) = j, X^{\pi}_{k}(\tau)\right] X^{\pi}_{k}(\tau) \mid j(t) = j \right] \nonumber \\ 
=& \sum^t_{\tau = 1} \sum_{k\in \KKK} \sum_{j\in \JJJ} p_j   \mathbb{E}\left[ \mathbf{1}\left( D_{ijk}\geq t - \tau + 1\right) A_{ijk} \right] x_{jk}(\tau) 
\leq  c_i \nonumber
\end{align}
Similarly, by taking expectation over each of the reward constraints, we have $\mathbb{E}[\sum^T_{t=1} \sum_{k\in \KKK} W_{ij(t)k}(t) X^{\pi}_{k}(t)]= \sum^T_{t=1} \sum_{j \in \JJJ} \sum_{k\in \KKK} p_j w_{ijk} x_{jk}(t) $ for each $i\in \III_r$. Hence, the claim about the objective value is shown, and the Lemma is proved. 
$\hfill\square$

\subsection{Proof of Lemma \ref{lem:benchmark2}}\label{pf:lemma_benchmark2}
It is evident that $T \lambda_* \leq T \cdot\text{opt}(\text{LP-E})$, since the solution $x$ defined as $x_{jk}(\tau)=y^*_{jk}$ for each $j\in \JJJ, k\in \KKK$ and $\tau\in \{1, \ldots, T\}$ is feasible to (LP-E), and the objective value of $x$ in (LP-E) is precisely equal to $\lambda^*$. Hence, we only need to show $T \cdot\text{opt}(\text{LP-E})- d_{\max} w_{\max} \leq T \lambda_*$. Since $D_{ijk} \in [0,d_{\max}]$ almost surely for all $i, j, k$, the resource constraints in (LP-E) can be equivalently written as $\sum^t_{\tau = \max(t-d_{\max},1)}  \sum_{j\in \JJJ}\sum_{k\in \KKK} p_j \mathbb{E}[A_{ijk}\mathbf{1}(D_{ijk} \geq t - \tau + 1) ] x_{jk}(\tau) \leq c_i ~  \forall i\in \III_c,~ t\in \{1, \ldots T\}$. Likewise, the resource constraints in (LP-S) can be equivalently written as $\sum_{j\in \JJJ}\sum_{k\in \KKK} \sum^{d_\text{max}}_{\tau = 1} p_j \mathbb{E}[A_{ijk}\mathbf{1}(D_{ijk} \geq \tau) ]  y_{jk} \leq c_i ~\forall i\in \III_c$. 

Thus, the dual of the expected LP (LP-E) can be express as:
\begin{subequations}
\begin{alignat}{2}
&\text{(LP-E-D0)}: ~\min\limits_{\boldsymbol{\alpha}, \boldsymbol{\beta}, \boldsymbol{\rho}} ~ \sum^T_{t=1} \left(\sum_{j \in \JJJ} \beta_{jt} + \sum_{i \in \III_c} c_i \alpha_{it}\right) & \nonumber\\
\text{s.t.}  &~\beta_{jt} + p_j \sum_{i \in \III_c}  \sum^{\min\{t+d_{\max} -1,T\}}_{\tau=t} \mathbb{E}[A_{ijk} \mathbf{1}(D_{ijk} \geq  \tau -t + 1)] \alpha_{i \tau} - p_j\sum_{i \in \III_r} w_{ijk} \rho_i \geq 0   &\quad &\forall j \in \JJJ, ~k \in \KKK, \nonumber \\
&&&~t \in \{1, \ldots, T\}   \nonumber\\
&T\sum_{i \in \III_r} \rho_i \geq 1   & \nonumber\\
&\rho_i \geq 0  & \quad &\forall i \in \III_r \nonumber\\
&\alpha_{it} \geq 0   &\quad &\forall i \in \III_c, ~t \in \{1, \ldots, T\} \nonumber\\
&\beta_{jt} \geq 0      &\quad &\forall j\in \JJJ, ~t \in \{1, \ldots, T\} \nonumber.
\end{alignat}
\end{subequations}
which has the same optimal value as
\begin{subequations}
\begin{alignat}{2}
&\text{(LP-E-D)}: ~\min\limits_{\boldsymbol{\alpha}, \boldsymbol{\beta}, \boldsymbol{\rho}} ~ \sum^T_{t=1} \left(\sum_{j \in \JJJ} p_j\beta_{jt} + \sum_{i \in \III_c} c_i \alpha_{it}\right) & \nonumber\\
\text{s.t.}  &~\beta_{jt} + \sum_{i \in \III_c}  \sum^{\min\{t+d_{\max} -1,T\}}_{\tau=t} \mathbb{E}[A_{ijk} \mathbf{1}(D_{ijk} \geq  \tau -t + 1)] \alpha_{i \tau} - \sum_{i \in \III_r} w_{ijk} \rho_i \geq 0   &\quad &\forall j \in \JJJ, ~k \in \KKK, \nonumber \\
&&&~t \in \{1, \ldots, T\}   \nonumber\\
&T\sum_{i \in \III_r} \rho_i \geq 1   & \nonumber\\
&\rho_i \geq 0  & \quad &\forall i \in \III_r \nonumber\\
&\alpha_{it} \geq 0   &\quad &\forall i \in \III_c, ~t \in \{1, \ldots, T\} \nonumber\\
&\beta_{jt} \geq 0      &\quad &\forall j\in \JJJ, ~t \in \{1, \ldots, T\} \nonumber.
\end{alignat}
\end{subequations}

The dual of (LP-S) has the same optimal value as:
\begin{subequations}
\begin{alignat}{2}
\text{(LP-S-D)}: &~\min\limits_{\boldsymbol{\alpha}, \boldsymbol{\beta}, \boldsymbol{\rho}} ~ \sum_{j \in \JJJ} p_j \beta_j + \sum_{i \in \III_c} c_i \alpha_i & \nonumber \\
\text{s.t.}  &~\beta_j + \sum_{i \in \III_c} \sum^{d_{\max}}_{\tau = 1} \mathbb{E}[A_{ijk} \mathbf{1}(D_{ijk} \geq  \tau)] \alpha_i - \sum_{i \in \III_r} w_{ijk} \rho_i \geq 0   &\quad &\forall j \in \JJJ, ~k \in \KKK   \nonumber\\
&\sum_{i \in \III_r} \rho_i \geq 1   & \nonumber\\
&\rho_i \geq 0  &\quad &\forall i \in \III_r \nonumber\\
&\alpha_i \geq 0   &\quad &\forall i \in \III_c \nonumber\\
&\beta_j \geq 0      &\quad &\forall j\in \JJJ \nonumber.
\end{alignat}
\end{subequations}
To show the inequality $T \cdot\text{opt}(\text{LP-E})- d_{\max} w_{\max} \leq T \lambda_*$, we start by considering an optimal solution $(\boldsymbol{\alpha}^*, \boldsymbol{\beta}^*, \boldsymbol{\rho}^*)$ to (LP-S-D). We claim that the solution $(\boldsymbol{\bar{\alpha}} = (\bar{\alpha}_{i\tau})_{i\in \III_c, \tau\in \{1, \ldots, T\}}, \boldsymbol{\bar{\beta}} = (\bar{\beta}_{j t})_{j\in \JJJ, t\in \{1, \ldots, T\}}, \boldsymbol{\bar{\rho}} = (\bar{\rho}_i)_{i\in \III_r})$ defined as 
\begin{align*}
\bar{\alpha}_{i\tau} &= \alpha^*_{i} / T \qquad\qquad \text{for all }i \in \III_c, ~t \in \{1, \ldots, T\}\nonumber\\
\bar{\beta}_{jt} &= 
\begin{cases}
\beta^*_j / T& \quad \text{for all }j\in \JJJ, t\in \{1, \ldots, T - d_\text{max}\}, \\
w_\text{max}/T& \quad \text{for all } j\in \JJJ, t\in \{T-d_\text{max}+1, \ldots, T\},
\end{cases}\nonumber\\
\bar{\rho}_{i} &= \rho^*_i / T \qquad\qquad \text{for all }i \in \III_r
\end{align*}
is feasible to (LP-E-D). Now, note that the objective value of $(\boldsymbol{\bar{\alpha}}, \boldsymbol{\bar{\beta}}, \boldsymbol{\bar{\rho}})$ in (LP-E-D) is equal to 
$$
\frac{w_\text{max}d_{\text{max}}}{T} + \sum_{j\in \JJJ}\left(1 - \frac{d_{\text{max}}}{T}\right)\beta^*_j + \sum_{i\in \III_c}c_i\alpha^*_c \leq \frac{w_\text{max}d_{\text{max}}}{T} + \lambda^*,
$$
thus to establish the inequality $T \cdot\text{opt}(\text{LP-E})- d_{\max} w_{\max} \leq T \lambda_*$, it suffices to show the feasibility of $(\boldsymbol{\bar{\alpha}}, \boldsymbol{\bar{\beta}}, \boldsymbol{\bar{\rho}})$ to (LP-E-D). Firstly, by the optimality of $(\boldsymbol{\alpha}^*, \boldsymbol{\beta}^*, \boldsymbol{\rho}^*)$ to (LP-S-D), we know that $\sum_{i\in \III_r}\rho^*_i = 1$, thus it suffices to verify the feasibility to the first set of constraints in (LP-E-D). In the case of $t\in \{T-d_\text{max}+1, \ldots, T\}$, the constraint is feasible for every $j, k$, since $\bar{\beta}_{jt} = w_\text{max}/T$ while $\sum_{i\in \III_r} w_{ijk}\bar{\rho}_{i}\leq w_\text{max}/T$. In the case of $t\in \{1, \ldots, T - d_\text{max}\}$, the constraint is feasible for every $j, k$, since
\begin{subequations}
\begin{alignat}{2}
& \bar{\beta}_{jt} + \sum_{i \in \III} \sum^{\min\{t+d_{\max},T\}}_{\tau=t} \mathbb{E}[A_{ijk} \mathbf{1}(D_{ijk} \geq  \tau - t+1)] \bar{\alpha}_{i \tau} - \sum_{i \in \III} w_{ijk} \bar{\rho}_i \nonumber \\
= & \frac{1}{T}\left[\beta^*_{jt} + \sum_{i \in \III} \sum^{t+d_{\max}}_{\tau=t} \mathbb{E}[A_{ijk} \mathbf{1}(D_{ijk} \geq  \tau - t+1)] \alpha^*_{i \tau} - \sum_{i \in \III} w_{ijk} \rho^*_i \right] \nonumber \\
= & \frac{1}{T}\left[\beta^*_{jt} + \sum_{i \in \III} \sum^{d_{\max}}_{\tau=1} \mathbb{E}[A_{ijk} \mathbf{1}(D_{ijk} \geq  \tau )] \alpha^*_{i \tau} - \sum_{i \in \III} w_{ijk} \rho^*_i \right] \nonumber \\
= & \frac{1}{T}\left[\beta^*_{jt} + \sum_{i \in \III} \sum^{d_{\max}}_{t=1} \mathbb{E}[A_{ijk} t\mathbf{1}(D_{ijk} =t )] \alpha^*_{i \tau} - \sum_{i \in \III} w_{ijk} \rho^*_i \right] \nonumber \\
= & \frac{1}{T}\left[\beta^*_{jt} + \sum_{i \in \III} \mathbb{E}[A_{ijk} D_{ijk}] \alpha^*_{i \tau} - \sum_{i \in \III} w_{ijk} \rho^*_i \right] 
\geq  0. \nonumber
\end{alignat}
\end{subequations}
Altogether, the Lemma is proved.
$\hfill\square$

\subsection{Proof of Lemma \ref{lem:emp}}\label{pf:lemma_emp}
For cleanliness in notation, we omit the superscript $(q)$ throughout the proof, by abbreviating $\tau^{(q-1)}/2$ as $\tau$, $\hat{p}^{(q)}$ as $\hat{p}$, $\hat{\lambda}^{(q)}_*$ as $\hat{\lambda}_*$, and $\epsilon^{(q)}_A, \epsilon^{(q)}_B, \epsilon^{(q)}_C$ as $\epsilon_A, \epsilon_B, \epsilon_C$. We first show $\lambda_* \leq \hat{\lambda}_{*}  + \epsilon_A + \epsilon_B$. Consider the solution $$
\hat{y}_{jk} = \frac{1}{1 + \epsilon_A}y^*_{jk}\text{ for each $j\in \JJJ,k\in \KKK$.}
$$
We claim that, with probability at least \textcolor{black}{$1 - \delta$}, $\hat{y}$ is feasible to (LP-RS)$^{(q)}$ and satisfies
\begin{align}
\sum_{j \in \JJJ} \sum_{k \in \KKK} \hat{p}_{j} w_{ijk} \hat{y}_{jk} 
\geq & \lambda_* - \epsilon_B - \epsilon_A. \label{eq:emp_step_1}
\end{align}
These two properties clearly implies inequality $\lambda_* \leq \hat{\lambda}_{*}  + \epsilon_A + \epsilon_B$ holds with probability at least \textcolor{black}{$1 - 2\delta$}. Thus, we focus on establishing the feasibility and (\ref{eq:emp_step_1}).

To verify the feasibility, we consider a fixed $i\in \III_c$, and apply the multiplicative Chernoff inequality (\ref{eq:multi_chernoff_+}) with $X_t = \sum_{k\in \KKK}v_{i, j(t), k} y^*_{j(t), k}$ for $t\in \{\tau+1, \ldots,  2\tau \}$ and $B = \gamma$. Note that $\frac{1}{\tau} \sum^{2\tau}_{s=\tau+1} X_s = \sum_{j \in \JJJ} \sum_{k\in \KKK} \hat{p}_{j} v_{ijk}y^*_{jk}$ and $\mathbb{E}[\frac{1}{\tau} \sum^{2\tau}_{s=(\tau+1} X_s] = \sum_{j \in \JJJ} \sum_{k\in \KKK} p_{j} v_{ijk} y^*_{jk}$. The Chernoff inequality implies that, with probability at least \textcolor{black}{$ 1 - \frac{\delta}{|\III_c|}$}, we have
\begin{align}
&\sum_{j \in \JJJ} \sum_{k\in \KKK} \hat{p}_{j} v_{ijk}y^*_{jk}\nonumber\\
 \leq  & \sum_{j \in \JJJ} \sum_{k\in \KKK} p_{j} v_{ijk} y^*_{jk} +  \sqrt{\frac{4\gamma}{\tau}\sum_{j \in \JJJ} \sum_{k\in \KKK} p_{j} v_{ijk} y^*_{jk} \log\frac{|\III_c|}{\delta}} + \frac{2\gamma}{\tau} \log\frac{|\III_c|}{\delta}\nonumber\\
\leq & c_i + \sqrt{\frac{4\gamma c_i}{\tau}\log\frac{|\III_c|}{\delta}} + \frac{2\gamma}{\tau} \log\frac{|\III_c|}{\delta} \nonumber \\
= & \textcolor{black}{ c_i (1 + \epsilon_A)}, \nonumber
\end{align}
which implies that $\sum_{j \in \JJJ} \sum_{k\in \KKK} \hat{p}_{j} a_{ijk} d_{ijk}\hat{y}_{jk} \leq c_i$. The union bound thus shows that $\hat{y}$ is feasible to (LP-RS)$^{(q)}$ with probability at least $1 - \delta$.

Similarly, for a fixed $i\in \III_r$, by applying the Chernoff inequality with  $X_t = \sum_{k\in \KKK}w_{i, j(t), k} y^*_{j(t), k}$ for $t\in \{\tau+1, \ldots, 2\tau \}$ and $B = w_\text{max}$, with probability at least  \textcolor{black}{$ 1 - \frac{\delta}{|\III_c|}$} we have
\begin{align}
\sum_{j \in \JJJ} \sum_{k \in \KKK} \hat{p}_{j} w_{ijk} y^*_{jk} 
\geq & \lambda_* - w_{\max} \sqrt{\frac{2 \log(|\III_r|/\delta)}{\tau}} = \lambda_* - \epsilon_B. \nonumber
\end{align}
Thus, inequality (\ref{eq:emp_step_1}) follows from the definition of $\hat{y}$.

Next, we show \textcolor{black}{$\Pr(\hat{\lambda}_* -  \epsilon_C  \leq \lambda_*)\geq 1 - \delta$}. Now, the dual of (LP-S) has the same optimal value as
\begin{subequations}
\begin{alignat}{2}
\text{(LP-S-D) }: \min\limits_{\boldsymbol{\alpha}, \boldsymbol{\beta}, \boldsymbol{\rho}}  &~ \sum_{j \in \JJJ} p_j \beta_j + \sum_{i \in \III_c} c_i \alpha_i & \nonumber\\
\text{s.t.}  &~\beta_j + \sum_{i \in \III_c} v_{ijk} \alpha_i - \sum_{i \in \III_r} w_{ijk} \rho_i \geq 0   &\quad &\forall j \in \JJJ, ~k \in \KKK   \nonumber\\
&\sum_{i \in \III_r} \rho_i \geq 1   & \nonumber\\
&\rho_i \geq 0   &\quad &\forall i \in \III_r \nonumber\\ 
&\alpha_i \geq 0   &\quad &\forall i \in \III_c \nonumber\\
&\beta_j \geq 0      &\quad &\forall j\in \JJJ \nonumber.
\end{alignat}
\end{subequations}
Like-wise, the dual of $\text{(LP-RS)}^{(q)}$ has the same optimal value as that of the following LP:
\begin{subequations}
\begin{alignat}{2}
\text{(LP-RS-D)}^{(q)}: ~\min\limits_{\hat{\boldsymbol{\alpha}}, \hat{\boldsymbol{\beta}}, \hat{\boldsymbol{\rho}}}  &~ \sum_{j \in \JJJ} \hat{p}_{j} \hat{\beta}_{j} + \sum_{i \in \III_c} c_i \hat{\alpha}_{i} & \nonumber\\
\text{s.t.}  &~\hat{\beta}_{j} + \sum_{i \in \III_c} v_{ijk} \hat{\alpha}_{i} - \sum_{i \in \III_r} w_{ijk} \hat{\rho}_{i} \geq 0   &\quad &\forall j \in \JJJ, ~k \in \KKK   \nonumber\\
&\sum_{i \in \III_r} \hat{\rho}_{i} \geq 1   & \nonumber\\
&\hat{\rho}_{i}  \geq 0   &\quad &\forall i \in \III_r \nonumber\\
&\hat{\alpha}_{i} \geq 0   &\quad &\forall i \in \III_c \nonumber\\
&\hat{\beta}_{j} \geq 0      &\quad &\forall j\in \JJJ \nonumber.
\end{alignat}
\end{subequations}
Since the feasible domains of $\text{(LP-S-D)}, \text{(LP-RS-D)}^{(q)}$ are the same, an optimal solution to (LP-S-D) is feasible to $\text{(LP-RS-D)}^{(q)}$. Consider a fixed optimal solutions $(\boldsymbol{\alpha}^{*}, \boldsymbol{\beta}^{*}, \boldsymbol{\rho}^{*})$ to (LP-S-D). On one hand, the objective value of $(\boldsymbol{\alpha}^{*}, \boldsymbol{\beta}^{*}, \boldsymbol{\rho}^{*})$ on (LP-S-D) is equal to $\lambda_* =\sum_{j \in \JJJ} p_j \beta^{*}_j + \sum_{i \in \III_c} c_i \alpha^{*}_i$. On the other hand, the objective value of $(\boldsymbol{\alpha}^{*}, \boldsymbol{\beta}^{*}, \boldsymbol{\rho}^{*})$ on $\text{(LP-RS-D)}^{(q)}$ can be used to bound $\hat{\lambda}_{*}$ as follows:
\begin{align}
\hat{\lambda}_{*} \leq & \sum_{j \in \JJJ} \hat{p}_{j} \beta^{*}_j + \sum_{i \in \III_c} c_i \alpha^{*}_i \nonumber \\
\leq & \sum_{j \in \JJJ} p_j \beta^{*}_j + \sum_{i \in \III_c} c_i \alpha^{*}_i + 2\sqrt{\frac{w_{\max} (\sum_{j \in \JJJ} p_j \beta^{*}_j)  \log(1/\delta)}{\tau} }+ 2\frac{w_{\max} \log(1/\delta)}{\tau} \quad \text{w.p. } 1-\delta \label{eq:emp_step_3}\\
\leq & \lambda_* + 2 w_{\max}\left[ \sqrt{\frac{ \log(1/\delta)}{\tau} }+ \frac{ \log(1/\delta)}{\tau}\right] \nonumber \\
= & \textcolor{black}{\lambda_* + \epsilon_C}. \nonumber 
\end{align}
Step (\ref{eq:emp_step_3}) is by the application of the Chernoff inequality in (\ref{eq:multi_chernoff_-}). The application crucially uses the fact that $\beta^*_j\leq w_\text{max}$, which is true because $\sum_{i\in \III_r}\rho^*_i = 1$ thanks to the optimality of the solution $(\boldsymbol{\alpha}^{*}, \boldsymbol{\beta}^{*}, \boldsymbol{\rho}^{*})$ to (LP-S-D). Altogether, the required bounds in the Lemma are proved. $\hfill\square$

\subsection{Proof of Lemma \ref{lem:iMWU}}\label{pf:lemma_iMWU}
The proof relies on a crucial application of Proposition \ref{prop:MWU}, with a judicious choice of $\ell(1), \ldots, \ell(\tau)$ (where we set $\tau = \tau^{(q-1)}$ ) that underpins the construction of Algorithm \ref{alg:IB}. Now, for each $s\in \{1, \ldots, \tau^{(q-1)}\}$, we define
\begin{equation}\label{eq:mwu_spec}
\ell_i(s)= \begin{cases}
w_{i j(s), k^{(q)}(s)} - (\hat{\lambda}^{(q)}_{*}- \epsilon^{(q)}_C) &\quad \text{ if $i\in \III_r$}\\
-v_{i, j(s), k^{(q)}(s)} + \text{min}\{c_i, v_\text{max}\}&\quad \text{ if $i\in \III_c$}.
\end{cases}
\end{equation}
It is evident that $|\ell_i(s)|\leq \gamma = \max\{w_\text{max}, v_\text{max}\}$ for all $i\in \III_r\cup \III_c, s\in \{1, \ldots, \tau^{(q-1)}/2\}$. In addition, under the specification of $\{\ell(s)\}^\tau_{s=1}$ in (\ref{eq:mwu_spec}), it can be directly verified that the MWU weigh vector $\vartheta(s)$ in (\ref{eq:vartheta}) (see Proposition \ref{prop:MWU}) is equal to  $({\bm \phi}^{(q)}(s),  {\bm \psi}^{(q)}(s))$ for each $s\in \{1, \ldots, \tau^{(q-1)}/2\}$. Applying Proposition \ref{prop:MWU} gives us the following inequalities (which simultaneously hold with certainty):
\begin{align}
\left[\frac{2}{\tau^{(q-1)}}\sum^{\tau^{(q-1)} / 2}_{s=1} w_{i, j(s), k^{(q)}(s)}\right] -\left( \hat{\lambda}^{(q)}_{*} - \epsilon^{(q)}_C\right) & \geq \Phi^{(q)} - \gamma\sqrt{\frac{8\log (|\III_r|+|\III_c|)}{\tau^{(q-1)}}}\quad\text{ for each $i\in \III_r$},\label{eq:mwu_Ir1}\\
-\left[\frac{2}{\tau^{(q-1)}}\sum^{\tau^{(q-1)} / 2}_{s=1} v_{i, j(s), k^{(q)}(s)}\right] + c_i&  \geq \Phi^{(q)} - \gamma\sqrt{\frac{8\log (|\III_r|+|\III_c|)}{\tau^{(q-1)}}} \quad\text{ for each $i\in \III_c$},\label{eq:mwu_Ic1}
\end{align}
where
\begin{align}\label{eq:Phi}
\Phi^{(q)}=
\frac{2}{\tau^{(q-1)}}\sum^{\tau^{(q-1)} / 2}_{s=1}&\left[\sum_{\iota\in \III_r}\phi^{(q)}_\iota(s)\left(w_{\iota j(s), k^{(q)}(s)} - (\hat{\lambda}^{(q)}_{*}- \epsilon^{(q)}_C)\right)\right.\nonumber\\
&\qquad +\left. \sum_{\iota\in \III_c}\psi^{(q)}_\iota(s) \left(-v_{\iota, j(s), k^{(q)}(s)} + \min\{c_\iota, v_\text{max}\}\right)\right].
\end{align}
To complete the proof of Lemma \ref{lem:iMWU}, \textcolor{black}{recall
\begin{align}
\epsilon^{(q)}_D = & 8\gamma\sqrt{\frac{\log ((|\III_r|+|\III_c|) / \delta)}{\tau^{(q-1)}}} \nonumber \\
\geq &\gamma\sqrt{\frac{8\log (|\III_r|+|\III_c|)}{\tau^{(q-1)}}} + 2\gamma\sqrt{\frac{1}{\tau^{(q-1)}}\log\frac{1}{\delta}} + \gamma\sqrt{\frac{8}{\tau^{(q-1)}}\log\frac{|\III_c|}{\delta}}. \nonumber
\end{align}} 
Then, it suffices to have the following 3 inequalities:
\begin{align}
\Pr\left(\Phi^{(q)} \geq -2\gamma\sqrt{\frac{1}{\tau^{(q-1)}}\log\frac{1}{\delta}}\right) &\geq 1 - 2\delta,\label{eq:MWU_conc_1}\\
\Pr\left(\sum^{\tau^{(q-1)} / 2}_{s=1} w_{i, j(s), k^{(q)}(s)} \leq   \sum^{\tau^{(q-1)}/2}_{s=1} \sum_{j\in \JJJ}p_j w_{i,j,\kappa({\bm \phi}^{(q)}(s), {\bm \psi}^{(q)}(s), j)} + \gamma\sqrt{2\tau^{(q-1)}\log\frac{|\III_r|}{\delta}}\right) &\geq 1- \frac{\delta}{|\III_r|},\label{eq:MWU_conc_2}\\
\Pr\left(\sum^{\tau^{(q-1)} / 2}_{s=1} v_{i, j(s), k^{(q)}(s)} \geq \sum^{\tau^{(q-1)}/2}_{s=1} \sum_{j\in \JJJ}p_j v_{i,j,\kappa({\bm \phi}^{(q)}(s), {\bm \psi}^{(q)}(s), j)} - \gamma\sqrt{2\tau^{(q-1)}\log\frac{|\III_c|}{\delta}} \right) &\geq 1- \frac{\delta}{|\III_c|},\label{eq:MWU_conc_3}.
\end{align}
We first show (\ref{eq:MWU_conc_1}). To begin, recall that $y^* = \{y^*_{jk}\}_{j\in \JJJ, k\in \KKK}$ is an optimal solution to (LP-S). 
\begin{align}
\Phi^{(q)}\geq &
\frac{2}{\tau^{(q-1)}}\sum^{\tau^{(q-1)} / 2}_{s=1}\left[\sum_{\iota\in \III_r}\phi^{(q)}_\iota(s)\left(\left\{\sum_{k\in \KKK}w_{\iota j(s), k}y^*_{j(s),k}\right\} - (\hat{\lambda}^{(q)}_{*}- \epsilon^{(q)}_C)\right)\right.\nonumber\\
&\qquad \qquad \qquad +\left. \sum_{\iota\in \III_c}\psi^{(q)}_\iota(s) \left(-\left\{\sum_{k\in \KKK}v_{\iota, j(s), k}y^*_{j(s),k}\right\} + \min\{c_\iota, v_\text{max}\}\right)\right]\label{eq:MWU_conc_1_1}\\
\geq&
\frac{2}{\tau^{(q-1)}}\sum^{\tau^{(q-1)} / 2}_{s=1}\left[\sum_{\iota\in \III_r}\phi^{(q)}_\iota(s)\left(\left\{\sum_{j\in \JJJ}\sum_{k\in \KKK}p_j w_{\iota j, k}y^*_{j(s),k}\right\} - (\hat{\lambda}^{(q)}_{*}- \epsilon^{(q)}_C)\right)\right.\nonumber\\
& +\left. \sum_{\iota\in \III_c}\psi^{(q)}_\iota(s) \left(-\left\{\sum_{j\in \JJJ}\sum_{k\in \KKK}p_jv_{\iota, j, k}y^*_{j(s),k}\right\} + \min\{c_\iota, v_\text{max}\}\right)\right] - \gamma \sqrt{\frac{4}{\tau^{(q-1)}}\log\frac{1}{\delta}} \text{ w.p. $\geq 1-\delta$}\label{eq:MWU_conc_1_2}\\
\geq&
\frac{2}{\tau^{(q-1)}}\sum^{\tau^{(q-1)} / 2}_{s=1}\left[\sum_{\iota\in \III_r}\phi^{(q)}_\iota(s)\left(\left\{\sum_{j\in \JJJ}\sum_{k\in \KKK}p_j w_{\iota j, k}y^*_{j(s),k}\right\} - \lambda^*\right)\right.\nonumber\\
&+\left. \sum_{\iota\in \III_c}\psi^{(q)}_\iota(s) \left(-\left\{\sum_{j\in \JJJ}\sum_{k\in \KKK}p_jv_{\iota, j, k}y^*_{j(s),k}\right\} + \min\{c_\iota, v_\text{max}\}\right)\right] - \gamma \sqrt{\frac{4}{\tau^{(q-1)}}\log\frac{1}{\delta}} \text{ w.p. $\geq 1-\delta$}\label{eq:MWU_conc_1_3}\\
\geq& - \gamma \sqrt{\frac{4}{\tau^{(q-1)}}\log\frac{1}{\delta}}\label{eq:MWU_conc_1_4}.  
\end{align}
Step (\ref{eq:MWU_conc_1_1}) is by the choice of $k^{(q)}(s)$ in Line \ref{alg:MWU_opt_oracle} in the MWU subroutine (Algorithm \ref{alg:IB}), and the inequality in step (\ref{eq:MWU_conc_1_1}) holds with certainty. Step (\ref{eq:MWU_conc_1_2}) is by an application of the Azuma Hoeffding inequality (see Proposition \ref{prop:multi_Chernoff}), with filtration $\{{\cal F}(s)\}^{\tau^{(q-1)}/2}_{s=1}$ defined as ${\cal F}(s) = \sigma(\{\hat{\lambda}^{(q)}_*\}\cup \{j(\tau)\}^s_{\tau=1})$. The inequality in (\ref{eq:MWU_conc_1_2}) holds with probability $\geq 1-\delta$. Step (\ref{eq:MWU_conc_1_3}) is by Lemma \ref{lem:emp}, and the inequality holds with probability at least $1-\delta$. Step (\ref{eq:MWU_conc_1_4}) is by the feasibility of $y^*$ to (LP-S), and also the fact that $\sum_{j\in \JJJ}\sum_{k\in \KKK}p_jv_{\iota, j, k}y^*_{j(s),k} \leq v_\text{max}$. The inequality (\ref{eq:MWU_conc_1_4})  holds with certainty. Altogether, inequality (\ref{eq:MWU_conc_1}) is shown.

Finally, inequalities (\ref{eq:MWU_conc_2}, \ref{eq:MWU_conc_3}) both follows from the Azuma-Hoeffding inequality (see Proposition~\ref{prop:hoeffding}). Inequality (\ref{eq:MWU_conc_2}) can be shown by considering  $\{w_{i, j(s), k^{(q)}(s)} - \sum_{j\in \JJJ}p_j w_{i,j,\kappa({\bm \phi}^{(q)}(s), {\bm \psi}^{(q)}(s), j)} \}^{\tau^{(q-1)}/2}_{s=1}$, which is a the martingale difference sequence with respect to the filtration $\{{\cal F}(s)\}^{\tau^{(q-1)}/2}_{s=1}$ defined as ${\cal F}(s) = \sigma(\{\hat{\lambda}^{(q)}_*\}\cup \{j(\tau)\}^s_{\tau=1})$. Crucially, in the conditional expectation $\mathbb{E}[w_{i, j(s), k^{(q)}(s)}  | {\cal F}(s-1)]$, the expectation is purely over the randomness in $j(s)$ (also note that $k^{(q)}(s)$ conditioned on $j(s) \cup {\cal F}(s-1)$  is deterministic), since the weight vector $({\bm \phi(s)}, {\bm \psi(s)})$ is ${\cal F}(s-1)$-measurable. Likewise, inequality (\ref{eq:MWU_conc_3}) can be shown by considering  $\{v_{i, j(s), k^{(q)}(s)} - \sum_{j\in \JJJ}p_j v_{i,j,\kappa({\bm \phi}^{(q)}(s), {\bm \psi}^{(q)}(s), j)} \}^{\tau^{(q-1)}/2}_{s=1}$, which is also a martingale difference sequence with respect to the filtration $\{{\cal F}(s)\}^{\tau^{(q-1)}/2}_{s=1}$. Altogether the Lemma is proved.

At this closure, we remark that it is vital to construct $\hat{\lambda}^{(q)}_*$ with a disjoint set of samples from Algorithm \ref{alg:IB}. Indeed, if we suppose the contrary, and let there be a case when the Virtual MWU and (LP-RS)$^{(q)}$ share a sample (say $j(s)$), then the probability distribution of $j(s)$ conditioned on ${\cal F}(s-1)$ needs not be $\textbf{p}$, since the random variable $j(s)$ is correlated with $\hat{\lambda}^{(q)}_*$. \hfill $\square$

\subsection{Proof of Lemma \ref{lem:bounding_reject}}\label{pf:lema_bounding_reject}
We start by observing that 
$$
\mathbb{E}[L(t) ~|~ {\cal H}(\tau^{(q-1)})] \geq 1- \sum_{i\in \III_c} \underbrace{\mathbb{E}\left[\mathbf{1}\left(\sum^{t-1}_{\tau=\max\{t-d_\text{max}, 1\}} \hat{A}_i(\tau)\mathbf{1}(\hat{D}_i(\tau) \geq t - \tau+1)> c_i - a_\text{max} \right)~\bigg |~ {\cal H}(\tau^{(q-1)})\right]}_{=M^{(q)}_i(t)}.
$$
The Lemma is proved by established these three steps. Firstly, we demonstrate that for any $t\in \{\tau^{(q-1)} + 1 + d_\text{max}, \ldots, \tau^{(q)}\}$, any $i\in \III_c$ and any fixed $\varepsilon > 0 $, the inequality
\begin{equation}\label{eq:reject_bd_1}
M^{(q)}_i(t) \leq \frac{1}{(1+\varepsilon)^{\frac{c_i}{a_\text{max}} - 1}}\cdot \exp\left[\frac{\varepsilon}{a_\text{max}}\frac{1}{1 + \bar{\epsilon}^{(q)}_D+ \eta}\cdot\frac{2}{\tau^{(q-1)}}\sum^{\tau^{(q-1)} / 2}_{s=1} v_{i, j(s), k^{(q)}(s)} \right]
\end{equation}
holds with certainty. Note that the right hand side of (\ref{eq:reject_bd_1}) is independent of $t$. Secondly, we demonstrate that the inequality 
\begin{equation}\label{eq:reject_bd_2}
\exp\left[\frac{\varepsilon}{a_\text{max}}\frac{1}{1 + \bar{\epsilon}^{(q)}_D+ \eta}\cdot\frac{2}{\tau^{(q-1)}}\sum^{\tau^{(q-1)} / 2}_{s=1} v_{i, j(s), k^{(q)}(s)} \right] \leq \exp\left[\frac{\varepsilon}{a_\text{max}}\frac{c_i + \epsilon^{(q)}_D}{1 + \bar{\epsilon}^{(q)}_D+ \eta} \right]
\end{equation}
holds for all $i\in \III_c$ with probability at least \textcolor{black}{$1-3\delta$}. Thirdly, by setting $\varepsilon = \frac{\eta}{1 + \bar{\epsilon}^{(q)}_D}$, we demonstrate the inequality
\begin{equation}\label{eq:reject_bd_3}
\frac{1}{(1+\varepsilon)^{\frac{c_i}{a_\text{max}} - 1}}\cdot\exp\left[\frac{\varepsilon}{a_\text{max}}\frac{c_i + \epsilon^{(q)}_D}{1 + \bar{\epsilon}^{(q)}_D+ \eta} \right] \leq \frac{2}{|\III_c|}\sqrt{\xi},
\end{equation}
which holds with certainty since inequality (\ref{eq:reject_bd_3}) only involves deterministic parameters. Combining the inequalities (\ref{eq:reject_bd_1}, \ref{eq:reject_bd_2}, \ref{eq:reject_bd_3}) shows the Lemma. In the remaining, we prove (\ref{eq:reject_bd_1}, \ref{eq:reject_bd_2}, \ref{eq:reject_bd_3}).

Inequality (\ref{eq:reject_bd_1}) is shown by the following string of calculations, where all equalities and inequalities hold almost surely:
\begin{align}
&M^{(q)}_i(t) \nonumber\\
= & \mathbb{E}\left[\mathbf{1}\left(\sum^{t-1}_{\tau=\max\{t-d_\text{max}, 1\}} \hat{A}_i(\tau)\mathbf{1}(\hat{D}_i(\tau) \geq t - \tau+1)> c_i - a_\text{max} \right)~\bigg |~ {\cal H}(\tau^{(q-1)})\right] \nonumber\\
= & \mathbb{E}\left[\mathbf{1}\left((1 + \varepsilon)^{\sum^{t-1}_{\tau=\max\{t-d_\text{max}, 1\}} \frac{\hat{A}_i(\tau)}{a_\text{max}}\mathbf{1}(\hat{D}_i(\tau) \geq t - \tau+1)}> (1 + \varepsilon)^{\frac{c_i}{a_\text{max}} - 1} \right)~\bigg |~ {\cal H}(\tau^{(q-1)})\right] \nonumber\\
\leq & \frac{1}{(1 + \varepsilon)^{\frac{c_i}{a_\text{max}} - 1}} \cdot \mathbb{E}\left[(1 + \varepsilon)^{\sum^{t-1}_{\tau=\max\{t-d_\text{max}, 1\}} \frac{\hat{A}_i(\tau)}{a_\text{max}}\mathbf{1}(\hat{D}_i(\tau) \geq t - \tau+1)}~\bigg |~ {\cal H}(\tau^{(q-1)})\right] \label{eq:reject_bd_1_step_1}\\
= & \frac{1}{(1 + \varepsilon)^{\frac{c_i}{a_\text{max}} - 1}} \cdot \prod^{t-1}_{\tau=\max\{t-d_\text{max}, 1\}}\mathbb{E}\left[(1 + \varepsilon)^{ \frac{\hat{A}_i(\tau)}{a_\text{max}}\mathbf{1}(\hat{D}_i(\tau) \geq t - \tau+1)}~\bigg |~ {\cal H}(\tau^{(q-1)})\right] \label{eq:reject_bd_1_step_2}\\
\leq & \frac{1}{(1 + \varepsilon)^{\frac{c_i}{a_\text{max}} - 1}} \cdot \prod^{t-1}_{\tau=\max\{t-d_\text{max}, 1\}}\left(1 + \varepsilon \cdot \mathbb{E}\left[\frac{\hat{A}_i(\tau)}{a_\text{max}}\mathbf{1}(\hat{D}_i(\tau) \geq t - \tau+1)~\bigg |~ {\cal H}(\tau^{(q-1)})\right]\right) \label{eq:reject_bd_1_step_3}\\
\leq & \frac{1}{(1 + \varepsilon)^{\frac{c_i}{a_\text{max}} - 1}} \cdot \prod^{t-1}_{\tau=\max\{t-d_\text{max}, 1\}}\left(1 + \varepsilon \cdot \mathbb{E}[B(\tau)] \mathbb{E}\left[\frac{\tilde{A}_i(\tau)}{a_\text{max}}\mathbf{1}(\tilde{D}_i(\tau) \geq t - \tau+1)~\bigg |~ {\cal H}(\tau^{(q-1)})\right]\right) \nonumber\\
\leq & \frac{1}{(1 + \varepsilon)^{\frac{c_i}{a_\text{max}} - 1}} \cdot \prod^{t-1}_{\tau=\max\{t-d_\text{max}, 1\}}\exp\left( \varepsilon \cdot \mathbb{E}[B(\tau)] \mathbb{E}\left[\frac{\tilde{A}_i(\tau)}{a_\text{max}}\mathbf{1}(\tilde{D}_i(\tau) \geq t - \tau+1)~\bigg |~ {\cal H}(\tau^{(q-1)})\right]\right) \label{eq:reject_bd_1_step_4}\\
= & \frac{1}{(1 + \varepsilon)^{\frac{c_i}{a_\text{max}} - 1}} \cdot \exp\left(  \frac{\varepsilon}{a_\text{max}} \cdot \frac{1}{1+\bar{\epsilon}^{(q)}_D + \eta} \sum^{t-1}_{\tau=\max\{t-d_\text{max}, 1\}} \mathbb{E}\left[\tilde{A}_i(\tau)\mathbf{1}(\tilde{D}_i(\tau) \geq t - \tau+1)~\bigg |~ {\cal H}(\tau^{(q-1)})\right]\right) \label{eq:reject_bd_1_step_5}\\
\leq & \frac{1}{(1 + \varepsilon)^{\frac{c_i}{a_\text{max}} - 1}} \cdot \exp\left(  \frac{\varepsilon}{a_\text{max}} \cdot \frac{1}{1+\bar{\epsilon}^{(q)}_D + \eta} \sum^{t}_{\tau=t-d_\text{max}+1} \sum^{d_{\max}}_{s=t-\tau+1} \mathbb{E}\left[\tilde{A}_i(\tau)\mathbf{1}(\tilde{D}_i(\tau) =s)~\bigg |~ {\cal H}(\tau^{(q-1)})\right]\right) \nonumber\\
= & \frac{1}{(1 + \varepsilon)^{\frac{c_i}{a_\text{max}} - 1}} \cdot \exp\left(  \frac{\varepsilon}{a_\text{max}} \cdot \frac{1}{1+\bar{\epsilon}^{(q)}_D + \eta} \sum^{d_{\max}}_{s=1} \sum^{t}_{\tau=t-s+1}  \mathbb{E}\left[\tilde{A}_i(\tau)\mathbf{1}(\tilde{D}_i(\tau) =s)~\bigg |~ {\cal H}(\tau^{(q-1)})\right]\right) \nonumber\\
= & \frac{1}{(1 + \varepsilon)^{\frac{c_i}{a_\text{max}} - 1}} \cdot \exp\left(  \frac{\varepsilon}{a_\text{max}} \cdot \frac{1}{1+\bar{\epsilon}^{(q)}_D + \eta} \sum^{d_{\max}}_{s=1} \sum^{t}_{\tau=t-s+1}  \mathbb{E}\left[\tilde{A}_i(t)\mathbf{1}(\tilde{D}_i(t) =s)~\bigg |~ {\cal H}(\tau^{(q-1)})\right]\right) \label{eq:reject_bd_1_step_6}\\
= & \frac{1}{(1 + \varepsilon)^{\frac{c_i}{a_\text{max}} - 1}} \cdot \exp\left(  \frac{\varepsilon}{a_\text{max}} \cdot \frac{1}{1+\bar{\epsilon}^{(q)}_D + \eta} \sum^{d_{\max}}_{s=1} \mathbb{E}\left[\tilde{A}_i(t)\cdot s \cdot \mathbf{1}(\tilde{D}_i(t) =s)~\bigg |~ {\cal H}(\tau^{(q-1)})\right]\right) \nonumber\\
= & \frac{1}{(1 + \varepsilon)^{\frac{c_i}{a_\text{max}} - 1}} \cdot \exp\left( \frac{\varepsilon}{a_\text{max}} \cdot \frac{1}{1+\bar{\epsilon}^{(q)}_D + \eta} \cdot \mathbb{E}\left[\tilde{A}_i(t)\tilde{D}_i(t) ~\bigg |~ {\cal H}(\tau^{(q-1)})\right]\right) \nonumber\\
=& \frac{1}{(1+\varepsilon)^{\frac{c_i}{a_\text{max}} - 1}}\cdot \exp\left[\frac{\varepsilon}{a_\text{max}}\frac{1}{1 + \bar{\epsilon}^{(q)}_D+ \eta}\cdot\frac{2}{\tau^{(q-1)}}\sum^{\tau^{(q-1)} / 2}_{s=1} v_{i, j(s), k^{(q)}(s)} \right]\label{eq:reject_bd_1_step_7}.
\end{align}
Step (\ref{eq:reject_bd_1_step_1}) is by the Markov inequality. Step (\ref{eq:reject_bd_1_step_2}) follows from the joint independence of $\{\frac{\hat{A}_i(\tau)}{a_\text{max}}\mathbf{1}(\hat{D}_i(\tau) \geq t - \tau+1)\}^{t-1}_{\tau=\max\{t-d_\text{max}, 1\}}$ conditioned on ${\cal H}(\tau^{(q-1)})$, thanks to the coupling argument. Step (\ref{eq:reject_bd_1_step_3}) is by the fact that $(1+\varepsilon)^a \leq 1+\varepsilon \cdot a$  for all $a\in [0, 1], \varepsilon > 0$. Step (\ref{eq:reject_bd_1_step_4}) is by the inequality $1 + \varepsilon \leq e^{\epsilon}$ which holds for all $\varepsilon > 0$. In step (\ref{eq:reject_bd_1_step_5}), it is crucial to that the range of summation, namely $\{\max\{t-d_\text{max}, 1\}, \ldots, t-1\}$, lies inside the time interval of phase $q$, by our assumption that $t\in \{\tau^{(q-1)} + 1 + d_\text{max}, \ldots, \tau^{(q)}\}$. Now, recall Algorithm \ref{alg:unknown_IB} Line \ref{alg:tilde} that constructs $\tilde{k}(t)$ and the definition of $(\tilde{W}(t), \tilde{A}(t), \tilde{D}(t)) \sim {\cal O}_{j(t), \tilde{k}(t)}$ in our coupling argument. These two facts imply that $\{(\tilde{A}_i(\tau), \tilde{D}_i(\tau))\}^t_{\tau = \max\{t-d_\text{max}, 1\}}$ are iid conditioned on ${\cal H}(\tau^{(q-1)})$, which leads to step (\ref{eq:reject_bd_1_step_6}). Finally, step (\ref{eq:reject_bd_1_step_7}) follows from taking the conditional expectation and recalling the way $\tilde{k}(t)$ is constructed in Algorithm \ref{alg:unknown_IB}.

Inequality (\ref{eq:reject_bd_2}) follows from a direct application of Lemma \ref{lem:iMWU}. Inequality (\ref{eq:reject_bd_3}) follows by routine calculations and the assumptions on $\bar{\epsilon}^{(q)}_D, \eta$. \textcolor{black}{Now, we set $\varepsilon = \frac{\eta}{1 + \bar{\epsilon}^{(q)}_D}$. By Assumption \ref{ass:cap} we know $\eta \leq 1$, and therefore $\varepsilon \leq 1$.} Recall that $\bar{\epsilon}^{(q)}_D = \epsilon^{(q)}_D / c_\text{min}$, we have
\begin{align}
\frac{1}{(1+\varepsilon)^{\frac{c_i}{a_\text{max}} - 1}}\cdot\exp\left[\frac{\varepsilon}{a_\text{max}}\frac{c_i + \epsilon^{(q)}_D}{1 + \bar{\epsilon}^{(q)}_D+ \eta} \right] &\leq  (1+\varepsilon)\cdot \frac{1}{(1+\varepsilon)^{\frac{c_i}{a_\text{max}}}}\cdot \exp(\frac{\varepsilon}{1+\varepsilon}\cdot\frac{c_i}{a_\text{max}})\\
& = (1+\varepsilon)\left[\frac{e^\varepsilon}{(1+\varepsilon)^{1+\varepsilon}}\right]^{\frac{c_i}{a_\text{max}} \cdot \frac{1}{1+\varepsilon}} \nonumber\\
&\leq (1+\varepsilon)\cdot \exp\left[-\frac{\varepsilon^2}{(1+\varepsilon)(2+\varepsilon)} \cdot \frac{c_i}{a_\text{max}}  \right]\nonumber\\
&\leq (1+\varepsilon)\cdot \exp\left[-\frac{\varepsilon^2}{(1+\varepsilon)(2+\varepsilon)} \cdot \frac{1}{\xi}  \right]\label{eq:reject_bd_3_1}.
\end{align}
Recall from (\ref{eq:param}) that
\textcolor{black}{
$$
\eta = \sqrt{\xi \log \frac{|\III_c|}{\xi}} ,\quad \bar{\epsilon}^{(q)}_D = \frac{\epsilon^{(q)}_D}{c_\text{min}}\text{, where }\epsilon^{(q)}_D = 8\gamma\sqrt{\frac{\log ((|\III_r|+|\III_c|) / \delta)}{\tau^{(q-1)}}}.
$$
}
\textcolor{black}{Recalling our assumption that $\tau^{(q-1)} \geq \frac{384\gamma^2}{c_\text{min}^2}\log\frac{|\III_r| + |\III_c|}{\delta}$, } 
(\ref{eq:reject_bd_3_1}) can be upper bounded as
\textcolor{black}{
$$
(\ref{eq:reject_bd_3_1}) \leq 2\exp\left( -\frac{\eta^2}{6\xi(1+\bar{\epsilon}^{(q)}_D)^2} \right) = 2\exp\left( -\frac{\tau^{(q-1)} c_{\min}^2 \log(|\III_c|/\xi)}{384\gamma^2 \log ((|\III_r|+|\III_c|) / \delta) } \right) \leq \frac{2}{|\III_c|}\sqrt{\xi}.  
$$
}
Altogether, the Lemma is proved. \hfill $\square$

\subsection{Proof of Lemma \ref{lem:correlated}}\label{pf:lemma_correlated}
We first remark that the Lemma does not follow directly from the conditional Chernoff inequality in Lemma \ref{lem:multi_AZ}. Indeed, the random variables $\{L(t)\hat{W}_i(t)\}^{\tau^{(q)}}_{t=\bar{\tau}^{(q-1)} +1}$ are correlated even when we condition on ${\cal H}(\tau^{(q-1)})$. Instead, we apply Lemma \ref{lem:multi_AZ} on suitable subsets of $\{L(t)\hat{W}_i(t)\}^{\tau^{(q)}}_{t=\bar{\tau}^{(q-1)} +1}$ that partition $\{L(t)\hat{W}_i(t)\}^{\tau^{(q)}}_{t=\bar{\tau}^{(q-1)} +1}$. To this end, we define $N = \lceil \frac{\tau^{(q)} - \bar{\tau}^{(q-1)}}{d_{\text{max}} + 1}\rceil$. For $\ell\in \{1, \ldots, d_\text{max}+1\}$ and $n\in \{1,\ldots, N\}$, we define the time index $$
t(n ; \ell) = \bar{\tau}^{(q-1)} + \ell + (n-1)\cdot(d_\text{max} + 1).
$$
Clearly, we have $\{t(n;\ell)\}_{n\in \{1, \ldots, N\}, \ell\in \{1,\ldots, d_\text{max}+1\}} \supseteq\{\bar{\tau}^{(q-1)} +1, \ldots, \tau^{(q)}\}$. Cruically, we observe that for any $\ell\in \{1, \ldots, d_\text{max}+1\}$, the random variables in the collection
\begin{equation*}
\Psi(\ell) = 
\{L(t(n;\ell))\cdot \hat{W}_i(t(n;\ell))\}^{N}_{n=1}
\end{equation*} 
are independent and identically distributed conditioned on ${\cal H}(\tau^{(q-1)})$. 

We first show the conditional independence. For every $t$, the random variable $L(t)\hat{W}_i(t)$ is $\sigma(\{(\hat{W}(\tau),\hat{A}(\tau), \hat{D}(\tau))\}^{t}_{\tau = t-d_\text{max}} )$-measurable. More precisely, by the definition of $L(t), \hat{W}_i(t)$, we know that there is a deterministic function $g_i$ such that $L(t)\hat{W}_i(t) = g_i(\{(\hat{W}(\tau),\hat{A}(\tau), \hat{D}(\tau))\}^{t}_{\tau = t-d_\text{max}} )$, where $g_i$ does not vary with $t$ and only depends on $i$. Since the time indexes in $\Psi(\ell)$ are at least $d_\text{max}+1$ time steps apart, we know that for any two distinct $n, n'\in \{1, \ldots, N\}$, the time indexes sets $\{t(n;\ell) - d_\text{max} , \ldots, t(n, \ell)\}$ and $\{t(n';\ell) - d_\text{max} , \ldots, t(n', \ell)\}$ are disjoint. By observing that 
$\{(\hat{W}(\tau),\hat{A}(\tau), \hat{D}(\tau))\}^{\tau^{(q)}}_{\tau = \bar{\tau}^{(q-1)} -d_\text{max}}$ are independent conditioned on ${\cal H}(\tau^{(q-1)})$, we know that the random variables in $\Psi(\ell)$ are independent conditioned on ${\cal H}(\tau^{(q-1)})$. 

The identically distributed part follows from the fact that, for any $t\in \{\bar{\tau}^{(q-1)} + 1, \ldots, \tau^{(q)}\}$, we know that $\{t-d_\text{max}, \ldots, t\} \subset \{\tau^{(q-1)}+1, \ldots, \tau^{(q)}\}$. In addition, by the coupling argument on the construction of $\hat{W}, \hat{A}, \hat{D}$, we know that $\{(\hat{W}(t), \hat{A}(t), \hat{D}(t))\}^{\tau^{(q)}}_{t=\tau^{(q-1)} + 1}$ are identically distributed conditioned on ${\cal H}(\tau^{(q-1)})$, since all the time indexes in $\{\tau^{(q-1)}+1, \ldots, \tau^{(q)}\}$ belong to phase $q$. Since $L(t)\hat{W}_i(t) = g_i(\{(\hat{W}(\tau),\hat{A}(\tau), \hat{D}(\tau))\}^{t}_{\tau = t-d_\text{max}} )$ and $g_i$ does not vary with $t$, we know that $\{L(t)\hat{W}_i(t)\}^{\tau^{(q)}}_{t=\bar{\tau}^{(q-1)} +1}$ are identically distributed conditioned on ${\cal H}(\tau^{(q-1)})$, which in particular implies that the random variables in $\Psi(\ell)$ are identically distributed conditioned on ${\cal H}(\tau^{(q-1)})$. 

After establishing the claim that the random variables in $\Psi(\ell)$ are conditionally iid for any $\ell$, we apply the conditional Chernoff inequality (Lemma \ref{lem:multi_AZ}) on the random variables in $\Psi(\ell)$, along with ${\cal F} = {\cal H}(\tau^{(q-1)})$ and $B=w_\text{max}$. Summing (\ref{eq:bound_hat_W_and_L}) over $t\in \{t(n;\ell)\}^N_{n=1}$ gives us that, with probability at least \textcolor{black}{$1-7\delta$}, it holds that 
\begin{align}
& \sum^N_{n=1}\mathbb{E}[\hat{W}_i(t(n;\ell)) L(t(n;\ell)) ~|~ {\cal H}(\tau^{(q-1)})  ] = \mu_{-}(\ell) \nonumber\\
\geq & \left(1 - \textcolor{black}{3\sqrt{\xi \log \frac{|\III_c|}{\xi}}} \right) N \lambda_* - N (\epsilon^{(q)}_A + \epsilon^{(q)}_B+\epsilon^{(q)}_C+ \textcolor{black}{(1+\lambda^*)} \epsilon^{(q)}_D )\nonumber
\end{align}
for all $\ell\in \{1, \ldots, d_{\text{max}}+1\}$.  
Observe that $\mu_{-}(\ell) \leq w_\text{max}N$ for all $\ell\in \{1, \ldots, d_{\text{max}}+1\}$ almost surely. The conditional Chernoff inequality shows us that with probability \textcolor{black}{$\geq 1 - 8\delta$}, we have 
\begin{equation}\label{eq:conc_ell}
\sum^N_{n=1}\hat{W}_i(t(n;\ell)) L(t(n;\ell)) \geq \mu_{-}(\ell) - w_\text{max}\sqrt{2N\log\frac{d_\text{max} + 1}{\delta}}
\end{equation}
for all $\ell\in \{1, \ldots, d_\text{max}+1\}$. 
Summing (\ref{eq:conc_ell}) over $\ell$ gives
\begin{align}
&\sum^{\tau^{(q)}}_{t = \bar{\tau}^{(q-1)}+1} \hat{W}_i(t) L(t) \nonumber\\
\geq & \left(1 -\textcolor{black}{3\sqrt{\xi \log \frac{|\III_c|}{\xi}}} \right)(\tau^{(q)} - \bar{\tau}^{(q-1)}) \lambda_*- (\tau^{(q)} - \bar{\tau}^{(q-1)})(\epsilon^{(q)}_A + \epsilon^{(q)}_B+\epsilon^{(q)}_C+ \textcolor{black}{(1+\lambda^*)}\epsilon^{(q)}_D ) \nonumber\\
&\qquad - w_\text{max}(d_\text{max}+1)  \sqrt{2N\log\frac{d_\text{max} + 1}{\delta}} \qquad\qquad \text{with probability \textcolor{black}{$\geq 1 - 8\delta$}}\nonumber\\
\geq & \left(1 -\textcolor{black}{3\sqrt{\xi \log \frac{|\III_c|}{\xi}}} \right)(\tau^{(q)} - \tau^{(q-1)}) \lambda_*- (\tau^{(q)} - \tau^{(q-1)})(\epsilon^{(q)}_A + \epsilon^{(q)}_B+\epsilon^{(q)}_C+ \textcolor{black}{(1+\lambda^*)}\epsilon^{(q)}_D ) \nonumber\\
&\qquad - w_\text{max}(d_\text{max}+1)  \sqrt{4\frac{\tau^{(q)} - \tau^{(q-1)}}{d_\text{max} + 1}\log\frac{d_\text{max} + 1}{\delta}} - 2w_\text{max} d_\text{max} \qquad\qquad \text{a.s.}\nonumber\\
= & \left(1 -\textcolor{black}{3\sqrt{\xi \log \frac{|\III_c|}{\xi}}} \right)(\tau^{(q)} - \tau^{(q-1)}) \lambda_*- \tilde{O}(\sqrt{\tau^{(q)} - \tau^{(q-1)}})\qquad \text{a.s.}.\label{eq:ell_bound_1}
\end{align}
Step (\ref{eq:ell_bound_1}) is by the fact that $\epsilon^{(q)}_A + \epsilon^{(q)}_B+\epsilon^{(q)}_C+ \textcolor{black}{(1+\lambda^*)}\epsilon^{(q)}_D = \tilde{O}(1/\sqrt{\tau^{(q)} - \tau^{(q-1)}})$. 
Altogether, the Lemma is proved. \hfill $\square$

\subsection{A Generalization of Lemma \ref{lem:multi_AZ} and Proof}\label{pf:lemma_multi_AZ}
We demonstrate a generalized version of Lemma \ref{lem:multi_AZ} that provides high probability upper and lower bounds to a sum of conditionally independent random variables in the following. 
\begin{lemma}[Multiplicative Chernoff inequality with both sides]\label{lem:multi_AZ_full}
Suppose random variables $\{X_t\}^N_{t=1}$ satisfy the following properties:
\begin{enumerate}
\item $X_1, \ldots, X_N$ are jointly independent conditional on a $\sigma$-algebra ${\cal F}$,
\item $\Pr(X_t \in [0, B]\text{ for all $t\in \{1, \ldots, N\}$}) = 1$ for some $B\in \mathbb{R}_{>0}$.
\item There exists real numbers $\delta_{-}, \delta_{+}\in [0, 1]$ and $\mu_{-}, \mu_{+} \in \mathbb{R}_{> 0}$ such that
\begin{equation}\label{eq:mu_plus_minus}
\Pr\left(\mathbb{E}\left[\sum^N_{t=1} X_t\mid {\cal F}\right] > \mu_{+}\right) \leq \delta_+, \qquad \Pr\left(\mathbb{E}\left[\sum^N_{t=1} X_t\mid {\cal F}\right] < \mu_{-}\right) \leq \delta_-.
\end{equation}
\end{enumerate} 
Then the following concentration inequalities hold for any fixed but arbitrary $\delta\in (0,1)$:
\begin{align}
\Pr\left(\sum^N_{t=1} X_t > \mu_+ + 2\sqrt{B\mu_+ \log\frac{1}{\delta}} + 2B \log\frac{1}{\delta} \right) &\leq \delta + \delta_+, \label{eq:multi_AZ_plus}\\
\Pr\left(\sum^N_{t=1} X_t < \mu_- - \sqrt{2B\mu_- \log\frac{1}{\delta}}\right) &\leq \delta + \delta_- .\label{eq:multi_AZ_minus} 
\end{align}
\end{lemma}

\begin{proof}{Proof of Lemma \ref{lem:multi_AZ_full}}
We first proof the inequality (\ref{eq:multi_AZ_plus}). To ease the discussion, we denote 
\begin{equation}\label{eq:epsilon_plus}
\epsilon_+ = 2\sqrt{\frac{B}{\mu_+} \log\frac{1}{\delta}} + \frac{2B}{\mu_+} \log\frac{1}{\delta} . 
\end{equation}
To begin, we have
\begin{align}
&\Pr\left(\sum^N_{t=1} X_t > (1 + \epsilon_+)\mu_+ \right)\nonumber\\
\leq &\Pr\left(\sum^N_{t=1} X_t > (1 + \epsilon_+)\mu_+ \text{, and }  \mathbb{E}\left[\sum^N_{t=1} X_t\mid {\cal F}\right] \leq \mu_{+}\right) + \delta_+, \label{eq:plus_step_1}
\end{align}
where step (\ref{eq:plus_step_1}) is by the assumption 3 in the statement of the Lemma. Next, 
\begin{align}
&\Pr\left(\sum^N_{t=1} X_t > (1 + \epsilon_+)\mu_+ \text{, and }  \mathbb{E}\left[\sum^N_{t=1} X_t\mid {\cal F}\right] \leq \mu_{+}\right)\nonumber\\
 = & \Pr\left((1 + \epsilon_+)^{\sum^N_{t=1} \frac{X_t}{ B}} > (1 + \epsilon_+)^{\frac{(1 + \epsilon_+)\mu_+}{B}}\text{, and }  \mathbb{E}\left[\sum^N_{t=1} X_t\mid {\cal F}\right] \leq \mu_{+}\right) \nonumber\\
 \leq & \frac{1}{(1 + \epsilon_+)^{\frac{(1 + \epsilon_+)\mu_+}{B}}} \cdot \mathbb{E}\left[(1 + \epsilon_+)^{\sum^N_{t=1} \frac{X_t}{ B}} \cdot \mathbf{1}\left( \mathbb{E}\left[\sum^N_{t=1} X_t\mid {\cal F}\right] \leq \mu_{+} \right)\right]\label{eq:plus_step_2} \\
 \leq & \frac{1}{(1 + \epsilon_+)^{\frac{(1 + \epsilon_+)\mu_+}{B}}} \cdot \mathbb{E}\left[\prod^N_{t=1} \left(1 + \epsilon_+ \frac{X_t}{ B}\right)  \cdot \mathbf{1}\left( \mathbb{E}\left[\sum^N_{t=1} X_t\mid {\cal F}\right] \leq \mu_{+} \right)\right]\label{eq:plus_step_3}.
\end{align}
Step (\ref{eq:plus_step_2}) is by the Markov inequality. More precisely, it is by taking expectation over the inequality
\begin{align*}
&(1 + \epsilon_+)^{\frac{(1 + \epsilon_+)\mu_+}{B}} \cdot \mathbf{1}\left((1 + \epsilon_+)^{\sum^N_{t=1} \frac{X_t}{ B}} > (1 + \epsilon_+)^{\frac{(1 + \epsilon_+)\mu_+}{B}}\text{, and }  \mathbb{E}\left[\sum^N_{t=1} X_t\mid {\cal F}\right] \leq \mu_{+}\right) \nonumber\\
\leq  & (1 + \epsilon_+)^{\sum^N_{t=1} \frac{X_t}{ B}} \cdot \mathbf{1}\left( \mathbb{E}\left[\sum^N_{t=1} X_t\mid {\cal F}\right] \leq \mu_{+} \right),
\end{align*}
which holds with certainty. Step (\ref{eq:plus_step_3}) is by the fact that $(1+\epsilon)^a \leq (1+ \epsilon \cdot a)$ for any $\epsilon > 0$ and $a\in (0, 1)$, and the fact that $X_t/B\in (0, 1)$ almost surely (by Assumption 2).

To continue with (\ref{eq:plus_step_3}), we have
\begin{align}
&\mathbb{E}\left[\prod^N_{t=1} \left(1 + \epsilon_+ \frac{X_t}{ B}\right)  \cdot \mathbf{1}\left( \mathbb{E}\left[\sum^N_{t=1} X_t\mid {\cal F}\right] \leq \mu_{+} \right)\right] \nonumber\\
&\mathbb{E}\left[\mathbb{E}\left[ \prod^N_{t=1} \left(1 + \epsilon_+ \frac{X_t}{ B}\right)  \cdot \mathbf{1}\left( \mathbb{E}\left[\sum^N_{t=1} X_t\mid {\cal F}\right] \leq \mu_{+} \right)\mid {\cal F} \right]\right] \nonumber\\
& = \mathbb{E}\left[\mathbf{1}\left( \mathbb{E}\left[\sum^N_{t=1} X_t \mid {\cal F}\right]\leq \mu_{+} \right)\cdot \prod^N_{t=1}\left(  1 + \mathbb{E}\left[ \epsilon_+ \frac{X_t}{ B} \mid {\cal F}\right]\right) \right] \label{eq:plus_step_4}\\
&\leq \mathbb{E}\left[\mathbf{1}\left( \mathbb{E}\left[\sum^N_{t=1} X_t \mid {\cal F}\right]\leq \mu_{+} \right)\cdot \exp\left( \epsilon_+ \sum^N_{t=1}\mathbb{E}\left[ \frac{X_t}{ B} \mid {\cal F}\right]\right) \right] \label{eq:plus_step_5}\\
&\leq \exp\left(\frac{\epsilon_+ \mu_+}{ B}\right). \label{eq:plus_step_6}
\end{align}
Step (\ref{eq:plus_step_4}) is by the independence of $X_1, \ldots, X_N$ conditioned on ${\cal F}$ in Assumption 1, as well as the fact that the indicator random variable $\mathbf{1}\left( \mathbb{E}\left[\sum^N_{t=1} X_t \mid {\cal F}\right]\leq \mu_{+} \right)$ is ${\cal F}$-measurable. Step (\ref{eq:plus_step_5}) is by the fact that $1+\epsilon \leq e^\epsilon$ for all $\epsilon \in \mathbb{R}$. 

After that, applying the bound  (\ref{eq:plus_step_6}) to (\ref{eq:plus_step_3}) gives
\begin{align}
& \Pr\left(\sum^N_{t=1} X_t > (1 + \epsilon_+)\mu_+ \text{, and }  \mathbb{E}\left[\sum^N_{t=1} X_t\mid {\cal F}\right] \leq \mu_{+}\right) \leq \left(\frac{e^{\epsilon_+}}{(1 + \epsilon_+)^{1 + \epsilon_+}}\right)^{\frac{\mu_+}{B}}\nonumber\\
&\leq \exp\left(-\frac{\epsilon_+^2}{2 + \epsilon_+} \cdot \frac{\mu_+}{B}\right)\label{eq:plus_step_7},
\end{align}
where step (\ref{eq:plus_step_7}) follows from the following technical calculations:
\begin{align*}
\log\left(\left(\frac{e^{\epsilon_+}}{(1+\epsilon_+)^{(1+\epsilon_+)}}\right)^{\frac{\mu_+}{B}}\right) &= \frac{\mu_+}{B} [\epsilon_+ - (1+\epsilon_+) \log(1+\epsilon_+)]\nonumber\\
&\leq \frac{\mu_+}{B} \left[\epsilon_+ - (1+\epsilon_+) \frac{\epsilon_+}{1+(\epsilon_+/2)}\right]\leq -\frac{\epsilon_+^2}{2 + \epsilon_+} \cdot \frac{\mu_+}{B}.
\end{align*}
To complete the proof of inequality (\ref{eq:multi_AZ_plus}), it remains to show that $\exp\left(-\frac{\epsilon_+^2}{2 + \epsilon_+} \cdot \frac{\mu_+}{B}\right)\leq \delta_+$, or equivalently $\frac{\epsilon_+^2}{2 + \epsilon_+} \cdot \frac{\mu_+}{B}\geq \log \frac{1}{\delta_+}. $ The inequality is evident from the definition of $\epsilon_+$. In the case of $\epsilon_+ \leq 2$, we have
$$
\frac{\epsilon_+^2}{2 + \epsilon_+} \cdot \frac{\mu_+}{B} \geq \frac{\epsilon_+^2}{4} \cdot \frac{\mu_+}{B} \geq \frac{1}{4}\left[ \frac{1}{\mu_+}\left(2\sqrt{B\mu_+ \log\frac{1}{\delta}}\right)\right]^2\cdot \frac{\mu_+}{B}  = \log\frac{1}{\delta_+}.
$$
In the case of $\epsilon_+ > 2$, we have
$$
\frac{\epsilon_+^2}{2 + \epsilon_+} \cdot \frac{\mu_+}{B} \geq \frac{\epsilon_+^2}{2\epsilon_+} \cdot \frac{\mu_+}{B} \geq \frac{1}{2} \cdot \frac{2B}{\mu_+} \log \frac{1}{\delta_+} \cdot \frac{\mu_+}{B}  = \log\frac{1}{\delta_+},
$$
and altogether (\ref{eq:multi_AZ_plus}) is proved. 

The proof of (\ref{eq:multi_AZ_minus}) is analogous to the previous proof of (\ref{eq:multi_AZ_plus}), but we include the former for completeness sake. To ease the discussion, we denote 
\begin{equation}\label{eq:epsilon_minus}
\epsilon_- = \sqrt{\frac{2B}{\mu_-} \log\frac{1}{\delta_-}}. 
\end{equation}
If $\epsilon_- \geq 1$, then (\ref{eq:multi_AZ_minus}) is clearly true since then $ \mu_- - \sqrt{2B\mu_- \log\frac{1}{\delta}}\leq 0$. Thus, wlog we assume $\epsilon_-\in (0, 1)$ in the following calculations. To begin, we have
\begin{align}
&\Pr\left(\sum^N_{t=1} X_t < (1 - \epsilon_-)\mu_- \right)\nonumber\\
\leq &\Pr\left(\sum^N_{t=1} X_t < (1 - \epsilon_-)\mu_- \text{, and }  \mathbb{E}\left[\sum^N_{t=1} X_t\mid {\cal F}\right] \geq \mu_{-}\right) + \delta_-, \label{eq:minus_step_1}
\end{align}
where step (\ref{eq:minus_step_1}) is by the assumption 3 in the statement of the Lemma. Next, 
\begin{align}
&\Pr\left(\sum^N_{t=1} X_t < (1 - \epsilon_-)\mu_- \text{, and }  \mathbb{E}\left[\sum^N_{t=1} X_t\mid {\cal F}\right] \geq \mu_{-}\right)\nonumber\\
 = & \Pr\left((1 - \epsilon_-)^{\sum^N_{t=1} \frac{X_t}{ B}} > (1 - \epsilon_-)^{\frac{(1 - \epsilon_-)\mu_-}{B}}\text{, and }  \mathbb{E}\left[\sum^N_{t=1} X_t\mid {\cal F}\right] \geq \mu_{-}\right) \nonumber\\
 \leq & \frac{1}{(1 - \epsilon_-)^{\frac{(1 - \epsilon_-)\mu_-}{B}}} \cdot \mathbb{E}\left[(1 - \epsilon_-)^{\sum^N_{t=1} \frac{X_t}{ B}} \cdot \mathbf{1}\left( \mathbb{E}\left[\sum^N_{t=1} X_t\mid {\cal F}\right] \geq \mu_{-} \right)\right]\label{eq:minus_step_2} \\
 \leq & \frac{1}{(1 - \epsilon_-)^{\frac{(1 - \epsilon_-)\mu_-}{B}}} \cdot \mathbb{E}\left[\prod^N_{t=1} \left(1 - \epsilon_- \frac{X_t}{ B}\right)  \cdot \mathbf{1}\left( \mathbb{E}\left[\sum^N_{t=1} X_t\mid {\cal F}\right] \geq \mu_{-} \right)\right]\label{eq:minus_step_3}.
\end{align}
Step (\ref{eq:minus_step_2}) is by the Markov inequality, similar to that in (\ref{eq:plus_step_2}). Step (\ref{eq:minus_step_3}) is by the fact that $(1-\epsilon)^a \leq (1- \epsilon \cdot a)$ for any $\epsilon , a\in (0, 1)$, and the fact that $X_t/B\in (0, 1)$ almost surely (by Assumption 2).

To continue with (\ref{eq:minus_step_3}), we have
\begin{align}
&\mathbb{E}\left[\prod^N_{t=1} \left(1 - \epsilon_- \frac{X_t}{ B}\right)  \cdot \mathbf{1}\left( \mathbb{E}\left[\sum^N_{t=1} X_t\mid {\cal F}\right] \geq \mu_{-} \right)\right] \nonumber\\
&\mathbb{E}\left[\mathbb{E}\left[ \prod^N_{t=1} \left(1 - \epsilon_- \frac{X_t}{ B}\right)  \cdot \mathbf{1}\left( \mathbb{E}\left[\sum^N_{t=1} X_t\mid {\cal F}\right] \geq \mu_{-} \right)\mid {\cal F} \right]\right] \nonumber\\
& = \mathbb{E}\left[\mathbf{1}\left( \mathbb{E}\left[\sum^N_{t=1} X_t \mid {\cal F}\right]\geq \mu_- \right)\cdot \prod^N_{t=1}\left(  1 - \mathbb{E}\left[ \epsilon_- \frac{X_t}{ B} \mid {\cal F}\right]\right) \right] \label{eq:minus_step_4}\\
&\leq \mathbb{E}\left[\mathbf{1}\left( \mathbb{E}\left[\sum^N_{t=1} X_t \mid {\cal F}\right]\leq \mu_{-} \right)\cdot \exp\left( -\epsilon_- \sum^N_{t=1}\mathbb{E}\left[ \frac{X_t}{ B} \mid {\cal F}\right]\right) \right] \label{eq:minus_step_5}\\
&\leq \exp\left(\frac{-\epsilon_- \mu_+}{ B}\right). \label{eq:minus_step_6}
\end{align}
Step (\ref{eq:minus_step_4}) is by the independence of $X_1, \ldots, X_N$ conditioned on ${\cal F}$ in Assumption 1, as well as the fact that the indicator random variable $\mathbf{1}\left( \mathbb{E}\left[\sum^N_{t=1} X_t \mid {\cal F}\right]\geq \mu_{-} \right)$ is ${\cal F}$-measurable. Step (\ref{eq:minus_step_5}) is by the fact that $1-\epsilon \leq e^{-\epsilon}$ for all $\epsilon \in \mathbb{R}$. 

After that, applying the bound  (\ref{eq:minus_step_6}) to (\ref{eq:minus_step_3}) gives
\begin{align}
& \Pr\left(\sum^N_{t=1} X_t < (1 - \epsilon_-)\mu_- \text{, and }  \mathbb{E}\left[\sum^N_{t=1} X_t\mid {\cal F}\right] \geq \mu_{-}\right) \leq \left(\frac{e^{-\epsilon_-}}{(1 - \epsilon_-)^{1 - \epsilon_-}}\right)^{\frac{\mu_-}{B}}\nonumber\\
&\leq \exp\left(-\frac{\epsilon_-^2}{2} \cdot \frac{\mu_-}{B}\right)\label{eq:minus_step_7},
\end{align}
where step (\ref{eq:minus_step_7}) follows from the following technical calculations:
\begin{align*}
\log\left(\left(\frac{e^{-\epsilon_-}}{(1-\epsilon_-)^{(1-\epsilon_-)}}\right)^{\frac{\mu_-}{B}}\right) &= \frac{\mu_-}{B} [-\epsilon_- - (1-\epsilon_-) \log(1-\epsilon_-)]\nonumber\\
&\leq \frac{\mu_-}{B} \left[- \epsilon_- +\left(\epsilon_- - \frac{\epsilon_-^2}{2}\right)\right]= -\frac{\epsilon_-^2}{2} \cdot \frac{\mu_-}{B}.
\end{align*}
To complete the proof of inequality (\ref{eq:multi_AZ_minus}), it remains to observe that $\exp\left(-\frac{\epsilon_-^2}{2} \cdot \frac{\mu_+}{B}\right)\leq \delta_-$, which in fact holds with equality by the definition of $\epsilon_-$. Altogether (\ref{eq:multi_AZ_minus}) is proved.  \hfill$\square$
\end{proof}

\newpage
\subsection{Table of notation} \label{tab:notation}

\begin{table}[h]
\centering
\resizebox{\textwidth}{!}{%
\begin{tabular}{|l|l|}
\hline
Notation                & Usage                                                                                                            \\ \hline
$w_{ijk}$               & mean reward under $i \in \III_r$, $j \in \JJJ$ and $k \in \KKK$; $w_{\max} = \max_{i,j,k} \{w_{ijk}\}$               \\ \hline
$v_{ijk}$                & mean resource volume $\mathbb{E}[A_{ijk}D_{ijk}]$ under $i \in \III_c$, $j \in \JJJ$, $k \in \KKK$; $v_{\max} = \max_{i,j,k} \{v_{ijk}\}$        \\ \hline
$\bar{\epsilon}^{(q)}_D$ & Error parameter describing loss caused by model uncertainty, used to slim the customer flow                \\ \hline
$\eta$                   & Discount parameter describing loss caused by parameter stochasticity, used to slim the customer flow       \\ \hline
$\hat{\lambda}^{(q)}_*$ & Estimate to benchmark $\lambda_*$ based on observations in time steps $(\tau^{(q-1)} / 2) + 1, \ldots, \tau^{(q-1)}$       \\ \hline
$\Theta^{(q)}$          & Weight vector set generated by iMWU over observations in time steps $1, \ldots, \tau^{(q-1)}/2$                  \\ \hline
$\eta(s)$               & Leaning rate of iMWU at time $s$                                                                                 \\ \hline
$\epsilon^{(q)}_C$       & Error parameter describing loss caused by model uncertainty, $\hat{\lambda}^{(q)}_{*} - \epsilon^{(q)}_C$ used as a lower bound for $\lambda_*$ \\ \hline
$\phi^{(q)}_{i}(s)$     & Weight vector regarding reward $i \in \III_r$, generated by iMWU in phase $q$ time step $s$                      \\ \hline
$\psi^{(q)}_{i}(s)$     & Weight vector regarding resource and usage duration $i \in \III_c$, generated by iMWU in phase $q$ time step $s$ \\ \hline
$\Gamma^{(q)}_i(s)$     & Exponent parameter regarding reward $i \in \III_r$, used to update $\phi^{(q)}_{i}(s)$                                              \\ \hline
$\Xi^{(q)}_i(s)$     & Exponent parameter regarding resource $i \in \III_c$, used to update $\psi^{(q)}_{i}(s)$                                              \\ \hline
\end{tabular}%
}
\end{table}

\end{APPENDICES}

\end{document}